\documentclass[onefignum,onetabnum]{siamart171218}

\usepackage{amsfonts}
\usepackage{graphicx}
\usepackage{booktabs}
\usepackage{epstopdf}
\usepackage{algorithmic}
\usepackage{diagbox}
\ifpdf
  \DeclareGraphicsExtensions{.eps,.pdf,.png,.jpg}
\else
  \DeclareGraphicsExtensions{.eps,.png,.jpg}
\fi

\newsiamremark{remark}{Remark}
\newsiamremark{hypothesis}{Hypothesis}
\crefname{hypothesis}{Hypothesis}{Hypotheses}
\newsiamthm{claim}{Claim}

\headers{ORAS for the CPM}{I. C. T. May, R. D. Haynes, and S. J. Ruuth}

\title{Schwarz solvers and preconditioners for the closest point method}

\author{Ian C. T. May\thanks{Simon Fraser University, Burnaby, British Columbia V5A 1S6, Canada (mayianm@sfu.ca)}
\and Ronald D. Haynes\thanks{Memorial University of Newfoundland, St.~John's, Newfoundland A1C 5S7, Canada (rhaynes@mun.ca)}
\and Steven J. Ruuth\thanks{Simon Fraser University, Burnaby, British Columbia V5A 1S6, Canada (sruuth@sfu.ca)}}

\usepackage{amsopn}

\ifpdf
\hypersetup{
  pdftitle={Schwarz solvers and preconditioners for the closest point method},
  pdfauthor={I. May, R. D. Haynes, and S. J. Ruuth}
}
\fi

\providecommand{\mat}[1]{\ensuremath { {\rm \bf #1} } }
\providecommand{\restrict}[2]{\ensuremath \left. #1 \right|_{#2}}
\providecommand{\surf}{\ensuremath \mathcal{S} }
\providecommand{\pdiff}[2]{\ensuremath \frac{\partial #1}{\partial #2} }
\DeclareMathOperator*{\argmin}{arg\,min}

\providecommand{\revd}[1]{ {#1} }
\providecommand{\revt}[1]{ {#1} }

\begin{document}

\maketitle

% REQUIRED
\begin{abstract}
  The discretization of surface intrinsic elliptic partial differential equations (PDEs) \revd{poses} interesting challenges not seen in flat space. The discretization of these PDEs typically proceeds by either parametrizing the surface, triangulating the surface, or embedding the surface in a higher dimensional flat space. The closest point method (CPM) is an embedding method that represents surfaces using a function that maps points in the embedding space to their closest points on the surface. In the CPM, this mapping also serves as an extension operator that brings surface intrinsic data onto the embedding space, allowing PDEs to be numerically approximated by standard methods in a narrow tubular neighborhood of the surface. We focus here on numerically approximating the positive Helmholtz equation, $\left(c-\Delta_\surf\right)u=f,~c\in\mathbb{R}^+$ by the CPM \revd{paired with finite differences}. This yields a large, sparse, and non-symmetric system to solve. Herein, we develop restricted additive Schwarz (RAS) and optimized restricted additive Schwarz (ORAS) solvers and preconditioners for this discrete system. In particular, we develop a general strategy for computing overlapping partitions of the computational domain, as well as defining the corresponding Dirichlet and Robin transmission conditions. We demonstrate that the convergence of the ORAS solvers and preconditioners can be improved by using a modified transmission condition where more than two overlapping subdomains meet. Numerical experiments are provided for a variety of analytical and triangulated surfaces. We find that ORAS solvers and preconditioners outperform their RAS counterparts, and that using domain decomposition as a preconditioner gives faster convergence over using it as a solver, \revd{as expected}. The methods exhibit good parallel scalability over the range of process counts tested.
\end{abstract}

% REQUIRED

% REQUIRED

% 65F08 Preconditioners for iterative methods
% 65F10 Iterative methods for linear systems
% 65N22 Solution of discretized equations
% 65N55 Multigrid methods; domain decomposition
% 65Y05 Parallel computation

% -------------------------------------------------------------
\section{Introduction}
\label{intro}
Elliptic PDEs \revd{can} lead to large coupled systems of equations upon discretization. Care is required in solving such systems to maintain reasonable memory requirements and to yield solutions in a timely manner. The focus here is on the efficient numerical solution of elliptic PDEs intrinsic to surfaces. As a representative model problem we consider the surface intrinsic positive Helmholtz equation
\begin{equation}
  \left(c-\Delta_\surf\right)u = f,
  \label{eq:shiftPoiss}
\end{equation}
where $\Delta_\surf$ denotes the Laplace-Beltrami operator associated with the surface $\surf\subset\mathbb{R}^d$, and $c\in\mathbb{R}^+$ is a positive constant.
Such systems also appear for certain time-discretized diffusion and reaction-diffusion phenomena intrinsic to surfaces \cite{CBM:RDonPC}.
This equation also arises when finding the eigenvalues of the Laplace-Beltrami operator \cite{CBM:Eig}, a problem that arises in shape classification \cite{Reuter:ShapeDNA} and other applications.

Several methods have been proposed for the discretization of surface intrinsic PDEs. Parametrizing the surface leads to efficient methods \cite{FloaterHormann:Para} when the parametrization is available or can be approximated. Finite element methods operating on a triangulation of the surface \cite{DziukElliot} lead to sparse and symmetric systems of equations, but the approach is sensitive to the quality of the triangulation. Embedding methods seek to place the surface and the discretization in a higher dimensional flat space, $\mathbb{R}^d$, where more familiar numerical methods \revd{are available}. The surface can be posed as a level set of a function in the embedding space where finite differences in $\mathbb{R}^d$ may be used \cite{Cheng:LSM}, though the representation of open surfaces and the imposition of artificial boundary conditions at the boundary of the computational domain is non-trivial.

The closest point method (CPM) \cite{SJR:CPM,CBM:ICPM} is an embedding method which extends values on the surface into the embedding space via a mapping from points in the embedding space to their closest points on the surface. The discretization of the surface intrinsic PDE can then be achieved by a standard flat-space method. The efficient solution of the system of equations produced by this method as applied to equation \eqref{eq:shiftPoiss} is the focus of this paper. The CPM allows for the discretization over a wide range of surfaces with the same implementation. Triangulated surfaces \cite{CBM:LSE}, surfaces of varying co-dimension \cite{SJR:CPM}, and moving surfaces \cite{Argy:MoveSurf} can all be treated by the CPM. Section \ref{bg:cpm} recalls the CPM and establishes the notation and conventions used herein.

Domain decomposition (DD) methods have proven to be a powerful approach for the efficient solution of systems of linear equations, especially those arising from elliptic PDEs \cite{QuarteroniValli,ToselliWidlund}. In these methods, the large problem of interest is decomposed into many smaller and computationally cheaper problems. These subproblems can be rooted in purely algebraic formulations or obtained from continuous subproblems that can then be discretized. The solutions of these subproblems are assembled together to approximate the solution of the targeted large system, the accuracy of which is then iteratively improved to achieve a given tolerance \cite{SmithBjorstadGropp,DoleanNataf}. Furthermore, this procedure makes an excellent preconditioner for Krylov methods such as GMRES \cite{SmithBjorstadGropp,DoleanNataf}.

In many DD schemes these smaller subproblems can be decoupled and solved in parallel. A prior paper from the authors of this article showed some preliminary results on the application of classical Schwarz type DD solvers and preconditioners to the CPM \cite{May:DD25}. \revd{The present work} focuses on additive Schwarz type solvers and preconditioners that allow for parallelization, and in particular, optimized Schwarz methods which use enhanced boundary, or transmission, conditions on the subproblems to accelerate convergence of the iteration. \revd{This work introduces a boundary alignment procedure (Section \ref{subd:itf}), the use of multiple Robin weights in the ORAS method (Section \ref{trans:robfo:crosspt}), and evaluates the utility of the developed methods as preconditioners for GMRES. These extensions over the previous work are accompanied by substantially more numerical tests. The result is a more robust, and general method.}

The CPM and domain decomposition schemes are each reviewed briefly in Section \ref{bg}. Section \ref{subd} considers the desirable properties of the disjoint subdomains, and how they may be created. Section \ref{subd:over} builds an overlapping decomposition from these disjoint subdomains and addresses some concerns that arise from the atypical geometry of the problem. Then, Section \ref{trans} completes the formulation of the solution scheme by placing appropriate boundary conditions on the now known subdomains. Section~\ref{results} evaluates the performance of the method by first sweeping over the introduced parameters, and then testing the parallel scalability of the method for a fixed set of parameters.
% -------------------------------------------------------------
\section{Background}
\label{bg}
The CPM and additive Schwarz type DD methods are presented here to highlight important considerations and establish notation.

\subsection{The closest point method}
\label{bg:cpm}
The CPM allows for the discretization of surface intrinsic differential operators using standard Cartesian methods by first extending the solution into a narrow region surrounding the surface.
% and lying in a higher dimensional embedding space.
By keeping the solution constant in the surface normal direction, the gradient, divergence, and Laplacian over the embedding space all recover their surface intrinsic values when restricted to the surface \cite{SJR:CPM}.
To make this discussion concrete, we formulate a discretization of equation \eqref{eq:shiftPoiss} by the CPM paired with finite differences over $\mathbb{R}^d$.

Consider a surface, $\surf$, embeddable in $\mathbb{R}^d$. Define the closest point function for $\surf$ as
\begin{eqnarray}
  CP_\surf:\mathbb{R}^d&\rightarrow&\surf \\
  x&\mapsto&\underset{y\in\surf}{\argmin}||x-y||_2, \nonumber
\end{eqnarray}
which identifies the closest point on the surface, $y$, for any point, $x$, in the embedding space. To avoid introducing discontinuities into $CP_\surf$, we shall assume that the computational domain, $\Omega_D$,  has been chosen to consist only of points within a distance $\kappa_\infty^{-1}$ of the surface $\surf$, where $\kappa_\infty$ is an upper bound on the curvatures of $\surf$ \cite{Chu:Vari}.

A finite difference discretization of a surface PDE by the CPM is obtained by specifying a finite difference stencil to approximate the flat-space differential operator
and a degree, $p$, for the interpolating polynomials forming the discrete extension operator. We place a uniform grid with spacing $h$ in a suitable neighborhood of the surface
%(with $\Sigma_A\subset \Omega_D$ denoting the set of nodes in the grid)
(this will be made precise later). Assuming sufficient smoothness of $f$ and $CP_\surf$ \cite{Maerz/Macdonald:cpfunctions},
values of $f$ on the surface may be extended off the surface  in the normal direction  to $\mathcal{O}\left(h^{p+1}\right)$  accuracy
by evaluating a tensor product barycentric Lagrangian interpolant \cite{Tref:Bary}.

\begin{figure}[hbtp]
  \centering
  \includegraphics[width=0.45\linewidth]{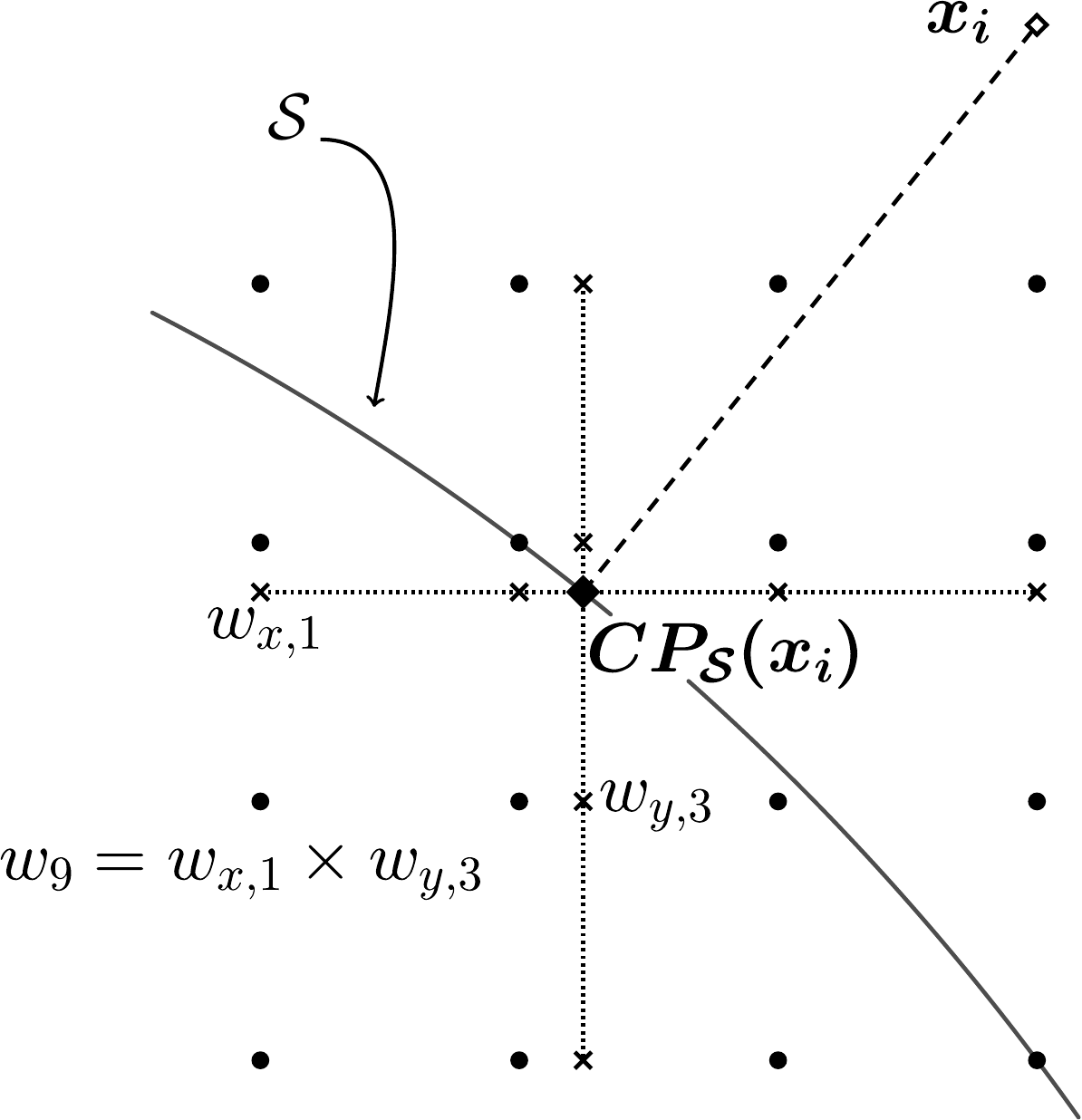}
  \caption{A bi-cubic extension stencil around a circle \revd{$\surf$} embedded in $d=2$ dimensional space.
Grid point $x_i$ (empty diamond) has closest point $CP_\surf(x_i)$ (filled diamond).  The associated 16 grid point interpolation stencil surrounding $CP_\surf(x_i)$ is displayed using filled circles.
The weights in the one dimensional stencils, shown by crosses and dotted lines, are given by barycentric Lagrangian interpolation.
The tensor product of the one dimensional interpolation weights gives the weights in the bi-cubic interpolation.
For instance, the $9^{\rm th}$ point in the extension stencil aligns with the first point in the $x-$direction and the third point in the $y-$direction, and thus has the weight $w_9=w_{x,1}\times w_{y,3}$.}
  \label{fig:extSten}
\end{figure}

To illustrate the procedure, Figure \ref{fig:extSten} shows an extension stencil in a $d=2$ dimensional embedding space using bi-cubic interpolation ($p=3)$.
In the figure, the $i^{\rm th}$ point in the grid, $x_i$,  is identified, along with its closest point, $CP_\surf(x_i)$, and the $(3+1)^2$ points surrounding it in the embedding space.
The $i^{\rm th}$ row of $\mat{E}$ will thus have 16 nonzero entries given by the tensor product of the weights in the one dimensional stencils used to approximate $f$ at $CP_\surf(x_i)$.

We are now in a position to define the grid.  The discrete extension operator relies only on points in a neighborhood
of the surface, and as a consequence the CPM may be posed over a neighborhood of the surface. This dramatically reduces the number of points in the mesh relative to a global computation (e.g., over a computational cube encompassing the shape), and ensures that during refinement the growth of the number of points will scale with the dimension of the surface rather than the dimension of the embedding space. In this work, we
%define a uniform grid in a neighborhood of the surface with spacing $h$ and
collect all nodes used for interpolation (the {\it active nodes}) into a set $\Sigma_A$ of size $N_A$. At the edge of this set there will be active nodes lacking neighboring nodes, and thus having incomplete finite difference stencils. We denote by $\Sigma_G$ the set of all {\it ghost nodes} (with $N_G=\#\Sigma_G$) needed to complete these stencils and note that the grid consists of the nodes in $\Sigma_A \cup \Sigma_G \subset \Omega_D$.

The sets $\Sigma_A$ and $\Sigma_G$ may be formed by assigning to them all nodes within specified radii of the surface \cite{SJR:CPM}.   We refer to this as the {\it tube-radius construction}.
Alternatively, the sets $\Sigma_A$ and $\Sigma_G$ may be formed algorithmically, as in \cite{CBM:ICPM}, where the active and ghost nodes are identified by the lists $L$ and $G$, respectively. This {\it algorithmic approach} leads to a smaller number of active nodes, and thus smaller linear systems. We use the former tube-radius construction for simple surfaces with analytically defined closest point functions, and the algorithmic approach for more complicated surfaces specified by an input triangulation. It is beneficial to demonstrate that both approaches yield domains appropriate for our proposed domain decomposition methods.

A possible discretization for the Laplace-Beltrami operator over $\surf$ may now be formed by composing the finite difference operator, $\Delta^h\in\mathbb{R}^{N_A \times(N_A+N_G)}$, with the discrete extension operator, $\mat{E}\in\mathbb{R}^{(N_A+N_G)\times N_A}$, containing the interpolation weights for all nodes in $\Sigma_A$ and $\Sigma_G$.
This {\it direct discretization} of the Laplace-Beltrami operator, $\Delta^h\mat{E}$, is a suitable  choice for explicit time-stepping of parabolic PDEs.
Unfortunately, the direct discretization has several small, spurious positive eigenvalues as well as complex eigenvalues with non-negligible imaginary parts. As a consequence, this discretization is not suitable for the numerical solution of eigenvalue problems \cite{CBM:Eig}. Furthermore, the spectrum of the direct discretization severely limits the allowable time step-size of implicit time-stepping methods, rendering them inefficient compared to the simpler explicit time-stepping methods \cite{CBM:ICPM}. Of final importance, use of the direct discretization tends to hinder the effectiveness of iterative solvers. One can resolve these issues by using a stabilized discretization
\begin{equation}
  \Delta_{\surf}^h = -\frac{2d}{h^2}I + \left(\Delta^h+\frac{2d}{h^2}I\right)\mat{E},
  \label{eq:impLapBel}
\end{equation}
which removes redundant extensions while remaining consistent \cite{CBM:ICPM,CBM:Eig}.

The effective iterative solution of linear systems defined by this discretization of the model problem~(\ref{eq:shiftPoiss}) sees two main obstacles: the matrix~(\ref{eq:impLapBel})  is non-symmetric even though the continuous operator is self-adjoint and, further, its spectrum has complex eigenvalues, though the imaginary parts are small.

% -------------------------------------------------------------
\subsection{Domain decomposition solvers and preconditioners}
\label{bg:dd}
There are many techniques for solving the linear systems that arise from the discretization of elliptic PDEs. If the system must be solved many times over and the operators involved are fixed (as occurs if using an implicit time-stepper with a fixed step size on a linear problem with non-varying coefficients), then prefactoring the matrix may be best (within time and memory constraints). Each time step then requires only a forward and backward solve with the initial cost of factorization being amortized by each step of time evolution. However, the solution of the same systems by adaptive time stepping methods, or the discretization of non-linear elliptic/parabolic systems solved by a Newton iteration, yield systems that are not amenable to prefactoring. Parallelization of iterative methods is much simpler than it is for direct solvers, and even the repeated solution of fixed systems may be better treated by iterative solvers with the use of appropriate preconditioners \cite{petsc-efficient,Saad:ItMeth}. As such, the design of effective iterative solvers is crucial and will be the focus herein.
% Similarly, the Jacobi-Davidson method for the eigenvalue problem involves solving an equation similar to \eqref{eq:shiftPoiss}, and is particularly attractive since this may be done through iterative methods \cite{blah:jacdav}, in contrast

A powerful option for solving such systems is through the use of DD methods as solvers, or as preconditioners for Krylov methods. Here the solution of the global problem is sought by repeatedly solving many smaller problems obtained by partitioning the domain and imposing appropriate conditions on the artificial boundaries. Thorough accounts of these methods can be found in \cite{ToselliWidlund,SmithBjorstadGropp,DoleanNataf,QuarteroniValli}, among many others. Here the focus will be on the restricted additive Schwarz method (RAS), and related methods as developed in \cite{Cai:RAS,Cyr:OMSORAS}.

The general formulation of these methods is best illustrated for the continuous problem \eqref{eq:shiftPoiss}, leaving the discretization for later. The global domain $\surf$ is divided into $N_S$ disjoint subdomains, $\widetilde{\surf}_j$, with a corresponding \revd{partitioning into $N_S$ overlapping subdomains $\surf_j$}. Define $\Gamma_{jk} := \partial\surf_j\cap\widetilde{\surf}_k$ as the portion of the boundary of the $j^{\rm th}$ overlapping \revd{subdomain} lying in the $k^{\rm th}$ disjoint \revd{subdomain}. Given an initial guess for the solution of the global problem $u^{(0)}$, defined at least over the boundaries of each subdomain $\surf_j$, a sequence of new solutions can be found by solving the subproblems
\begin{equation}
  \begin{cases}
    \left(c-\Delta_\surf\right)u_j = f, \quad&{\rm in}~\surf_j, \\
    \mathcal{T}_{jk}u_j = \mathcal{T}_{jk}u^{(n)}, \quad&{\rm on}~\Gamma_{jk},~~\forall~k,
  \end{cases}
  \label{eq:contORASIter}
\end{equation}
for the local solutions $u_j$.
A new global solution can be constructed from the restriction of the local solutions to their disjoint \revd{subdomains}.
\begin{equation}
  u^{(n+1)} = \sum\limits_j\restrict{u_j^{(n+1)}}{\widetilde{\surf}_j}.
  \label{eq:restCons}
\end{equation}
The boundary operators, $\mathcal{T}_{jk}$, transmit information between the subdomains, and are thus usually called transmission operators. The operators of interest herein are
\begin{eqnarray}
  \mathcal{T}_{jk} &=& identity\qquad{\rm or}\label{eq:dirTrans} \\
  \mathcal{T}_{jk} &=& \left(\pdiff{}{\hat{q}_{jk}} + \alpha\right),
  \label{eq:robTrans}
\end{eqnarray}
where \revd{$\hat{q}_{jk}$} is the outward pointing unit \textit{conormal} vector on $\Gamma_{jk}$. The conormal is defined as the vector that is simultaneously orthogonal to the surface normal, \revd{$\hat{n}$}, and the subdomain boundary, $\Gamma_{jk}$. The first option \eqref{eq:dirTrans} enforces Dirichlet conditions to produce the classic restricted additive Schwarz (RAS) solver \cite{Cai:RAS} when paired with the solution reconstruction in equation \eqref{eq:restCons}. The second option \eqref{eq:robTrans} enforces a Robin condition along each interface and gives the optimized restricted additive Schwarz family of solvers \cite{Cyr:OMSORAS}.

Following the notation for the continuous problems, we take $\Sigma_A$, $\widetilde{\Sigma}_j$, and $\Sigma_j$ as the global set, the disjoint subsets, and the overlapping subsets, respectively, of the active nodes in the discretization. Equation \eqref{eq:impLapBel} is gathered into the matrix $\mat{A}$ so that $\mat{A}\mat{u}=\mat{f}$ corresponds to a discretization of \eqref{eq:shiftPoiss} with $\mat{u}$ and $\mat{f}$ as vectors supported on $\Sigma_A$. The algebraic RAS iteration then proceeds by defining restriction operators $\widetilde{\mat{R}}_j:\Sigma_A\rightarrow\widetilde{\Sigma}_j$ and $\mat{R}_j:\Sigma_A\rightarrow\Sigma_j$, which simply truncate functions on $\Sigma_A$ to functions on $\widetilde{\Sigma}_j$ and $\Sigma_j$ respectively. Their transposes extend functions on these subdomains by zeros onto the entirety of $\Sigma_A$. Linearity of the problem may be exploited to instead solve for additive corrections to the current solution estimate \cite{DoleanNataf}, which naturally induces Dirichlet transmission conditions. For the RAS solver, the right hand side of the local problems becomes the restriction of the residual at the prior iteration,
\begin{equation}
  \mat{R}_j\mat{A}\mat{R}_j^T\mat{e}_j = \mat{R}_j\left(\mat{f}-\mat{A}\mat{u}^{(n)}\right)
  \label{eq:simpleRAS}
\end{equation}
where the local corrections, $\mat{e}_j$, are subsequently limited to the disjoint partition $\widetilde{\Sigma}_j$, \revd{extended by zeros onto $\Sigma_A$,} and added to the existing solution. Defining the local operator $\mat{A}_j = \mat{R}_j\mat{A}\mat{R}_j^T$, and gathering all steps together one finds the update
\begin{equation}
  \mat{u}^{(n+1)} = \mat{u}^{(n)} + \left[\sum\limits_{j=1}^{N_S}{\widetilde{\mat{R}}_j^T\mat{A}_j^{-1}\mat{R}_j}\right]\left(\mat{f}-\mat{A}\mat{u}^{(n)}\right),
  \label{eq:algRAS}
\end{equation}
where the summation in brackets forms the preconditioning operator $\mat{M}_{RAS}^{-1}$. This may be embedded within a Krylov scheme for additional acceleration, which effectively solves $\mat{M}_{RAS}^{-1}\mat{A}\mat{u}=\mat{M}_{RAS}^{-1}\mat{f}$ instead of the original system. Notably, the local solves may be done in parallel, requiring limited communication from the neighboring subdomains to populate the overlapping portion of the residuals. With the introduction of more advanced transmission conditions, the local operators $\mat{A}_j$ cannot be obtained from simple restrictions of the global matrix \cite{Cyr:OMSORAS,DoleanNataf} as in equation \eqref{eq:simpleRAS}, but the overall form of \eqref{eq:algRAS} as a solver or preconditioner remains valid.

% -------------------------------------------------------------
\section{Disjoint subdomain construction}
\label{subd}
To apply DD methods to the linear system arising from the CPM, an appropriate decomposition of the active nodes, $\Sigma_A$, must be found. The ability to obtain an arbitrary number of partitions for any given surface, relying only on information the CPM can provide, is critical to obtain a general method.

% -------------------------------------------------------------
\subsection{Graph based partitioning}
\label{subd:graph}
The  problem of disjoint subdomain construction is closely related to the graph partitioning problem and it is quite popular within the DD community to use graph partitioning tools such as METIS \cite{METIS} to partition the mesh.
The discrete domain, a point cloud in this case, is reinterpreted as a graph with connectivity describing the coupling between different nodes. For typical problems posed on flat space, this graph construction is fairly simple. For example,
finite difference schemes have nodes connected to their immediate neighbors, and possibly further connected through the stencil of the discrete operator. Finite element discretizations use elements as nodes and connect them to those elements sharing a face or edge (the dual-graph of the mesh \cite{METIS}).

The CPM offers several possible connectivity schemes due to the presence of extension, ambient discrete Laplacian, and assembled discrete Laplace-Beltrami operators.
Connecting nodes through their extension stencils (Figure \ref{fig:extSten}), or through the full stencil of the discrete Laplace-Beltrami operator, however, yields graphs with a large number of edges. The large number of edges present can create partitions that are contiguous with respect to the graph but not spatially.

Our partitioning strategy uses  the stencil of the ambient Laplacian.   This creates a graph which connects each node to its immediate neighbors (and no others). We have found that this creates a much simpler graph
than the other possible connectivity schemes.  Furthermore, it ensures that the generated partitions are contiguous within the graph as well as the mesh.

%-------------------------------------------------------------
\subsection{Alignment of the partition interfaces}
\label{subd:itf}
In the CPM, \revd{surface values are extended to the embedding space in a manner constant along the direction normal to the surface.} As a consequence, it is natural to seek partitions such that the interfaces between partitions are
aligned with the surface normals.  While the transmission conditions defined in Section \ref{trans} can be applied directly to the subdomains generated from the METIS splitting,
 the ignorance of the underlying geometry also means that the partitions will not respect the surface normals.

\begin{figure}
  \centering
  \includegraphics[width=0.95\linewidth]{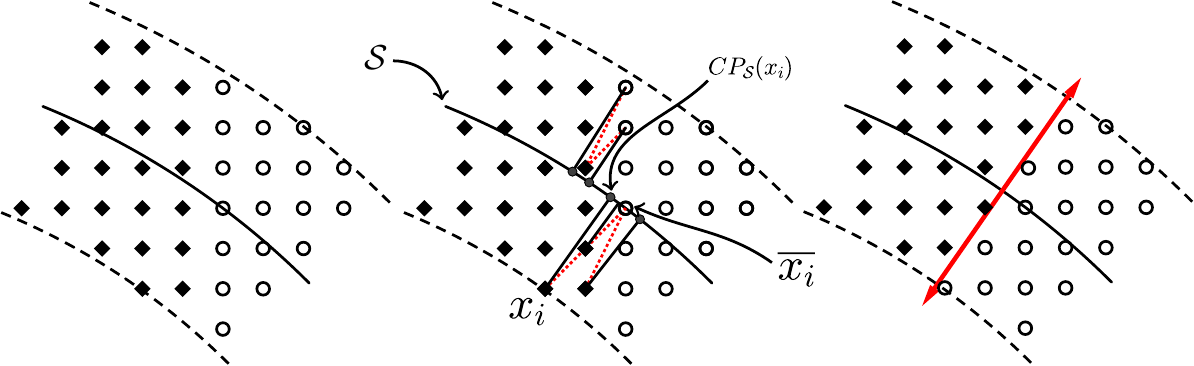}
  \caption{The left panel displays a typical disjoint partitioning from the nearest neighbor coupling scheme.  Diamonds denote one subdomain and circles the other. The middle panel displays select interface nodes and their closest points, and highlights the relationship with solid black lines. The corresponding nearest grid points are related to the interface nodes by red dashed lines. \revd{In the plot, we label the surface $\surf$}, an interface node $x_i$, its closest point $CP_\surf(x_i)$, and the corresponding nearest grid point $\overline{x}_i$. Since $\overline{x}_i$ belongs to a different subdomain than $x_i$, the interface node $x_i$ is marked for migration. The right panel illustrates the quality of the reassignment using a surface normal (red).}
  \label{fig:djtAlign}
\end{figure}

\begin{figure}[htpb]
  \centering
  \includegraphics[width=0.95\linewidth]{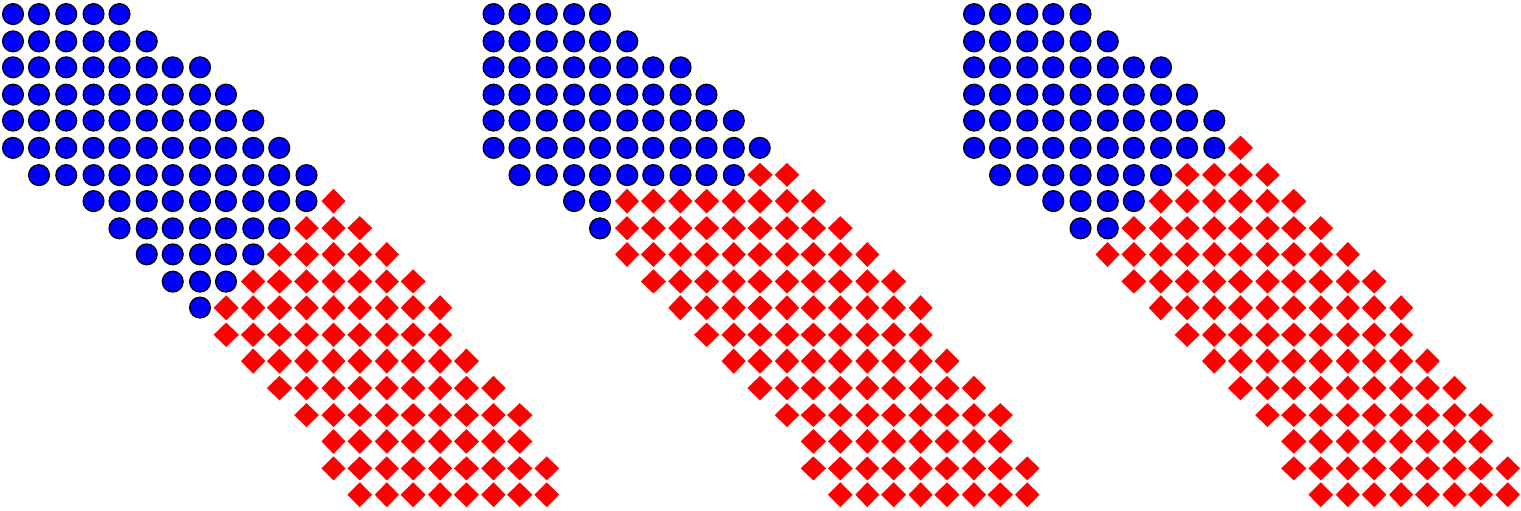}
  \caption{Partitions of the mesh for a circular problem. The left panel shows idealized partitions obtained  by manually splitting the nodes into wedges around the circle. The middle panel shows a pair of partitions obtained from the METIS splitting using the nearest neighbor coupling.  Notice that the interface does not respect the surface normal direction.   Applying  the process discussed in Section \ref{subd:itf} yields the result displayed in the right panel.  Improved alignment with the normal is obtained.}
  \label{fig:partExamples}
\end{figure}

To improve interface alignment, we propose a simple procedure; see Figure \ref{fig:djtAlign} for an illustration.
First, nodes along the interface between partitions are identified as those having neighbors in two or more partitions.
%For each interface node $x_i$, its closest point, $CP_\surf(x_i)$, is known.
For each interface node $x_i$, we determine the nearest grid point $\overline{x}_i$ to $CP_\surf(x_i)$.
(In the infrequent case where $CP_\surf(x_i)$ is a grid point, we assign $\overline{x}_i = CP_\surf(x_i)$.)
%will not coincide with a grid point, but can identified with the nearest grid point $\overline{x}_i$.
If the interface node $x_i$ belongs to a different partition than $\overline{x}_i$, then $x_i$ is flagged for migration into the partition containing $\overline{x}_i$. After all interface nodes are queried, those flagged are moved to their target partitions. This process may need to be repeated to fully align interfaces that are far from normal. We have found that applying $p+1$ passes ($p$ being the interpolation degree) works well in practice, and use this approach in all our numerical results. See Figure \ref{fig:partExamples} for a comparison of results using the METIS splitting and the proposed alignment procedure.

% -------------------------------------------------------------
\section{Overlapping subdomain construction}
\label{subd:over}
The procedure of Section~\ref{subd} partitions the global set of active nodes, $\Sigma_A$, into disjoint subsets, $\widetilde{\Sigma}_j, j=1\ldots,N_S$.
We have already seen in Section \ref{bg:dd} that overlapping partitions $\Sigma_j$ must be constructed from $\widetilde{\Sigma}_j$, just as the subdomains $\surf_j$ were constructed from $\widetilde{\surf}_j$ in the continuous formulation.

\revd{In Section \ref{trans} we will} provide full details on the construction of the \revd{discrete} local problems. Several sets of nodes will be required. These will be defined as they are needed in Sections~\ref{subd:bndNodes} and \ref{subd:bndGeo}.
Throughout, we shall refer to Figure~\ref{fig:sCons}, a partition for a circular arc, to help illustrate the main ideas.
Briefly, the sets to be constructed are:
\begin{itemize}
\item $\widetilde{\Sigma}_j, j=1\ldots,N_S$: disjoint sets of active nodes, which together form the global set $\Sigma_A$ (circles in Figure~\ref{fig:sCons}),
\item $\Sigma_j$: the overlapping set of active nodes in the local subproblem, $\widetilde{\Sigma}_j\subset\Sigma_j$ (circles and triangles in Figure \ref{fig:sCons}),
\item $\Sigma_j^G$: ghost nodes for the local subproblem that are also ghost nodes in the global problem (crossmarks in Figure~\ref{fig:sCons}),
\item $\Sigma_j^{BC}$: nodes needed to enforce the transmission condition (diamonds and plus signs in Figure \ref{fig:sCons}),
\item $\Lambda_j$: effective boundary locations, generated as closest points of the interface nodes in $\Sigma_j$ (approximately overlapping stars in Figure \ref{fig:sCons}).
\end{itemize}

% -------------------------------------------------------------
\subsection{Formation of overlaps and boundary nodes}
\label{subd:bndNodes}
\begin{figure}[htpb]
  \centering
  \includegraphics[width=0.6\linewidth]{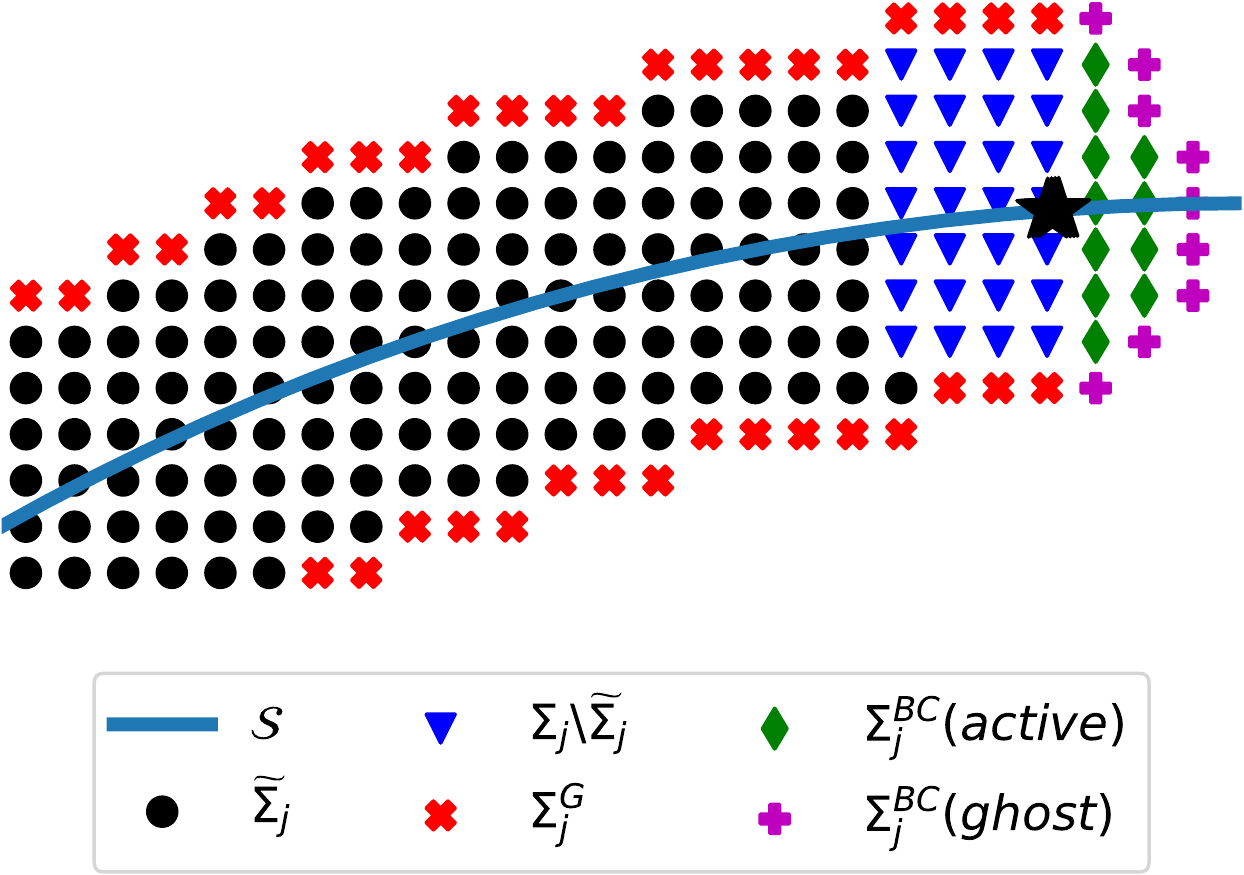}
  \caption{The node sets used in the construction of one subproblem.  These are built up successively. We start with the disjoint nodes $\widetilde{\Sigma_j}$ coming from the graph partitioning (black circles). Four layers of overlap nodes are added to create the grown partition $\Sigma_j$, which includes all of the black nodes as well as the blue triangles. The ghost node set $\Sigma_j^G$ is formed as a layer around these active nodes and is shown by red crossmarks. The nodes from incomplete stencils are shown as green diamonds, with their ghost nodes shown as purple plus signs, forming $\Sigma_j^{BC}$ collectively. Effective boundary locations, defined in Section \ref{subd:bndGeo}, are displayed as (approximately overlapping) black stars.}
  \label{fig:sCons}
\end{figure}

The set $\Sigma_j$ is grown from the disjoint partition $\widetilde{\Sigma}_j$ by layering additional, globally active nodes. The first layer is built by visiting all interface nodes in the local mesh $\Sigma_j$ and adding any missing neighboring nodes that are in $\Sigma_A$. The added nodes become the new interface nodes, and the process is repeated. A total of $N_O$ passes through the interface nodes creates the desired overlap between the local \revd{partitions}. An illustration of these overlap nodes (blue triangles) with $N_O=4$ is shown in Figure \ref{fig:sCons}.

The set of local ghost nodes, $\Sigma_j^G$, is formed from the nodes needed to complete the finite difference stencils over the active nodes $\Sigma_j$. The set $\Sigma_j^G$ (crossmarks in Figure~\ref{fig:sCons}) consists only of those nodes that are ghost nodes in the global problem as well.
%These are displayed as crossmarks in Figure~\ref{fig:sCons}.
Together, $\Sigma_j$ and $\Sigma_j^G$ define the interior of the $j^{\rm th}$ subdomain. Over these nodes, the discretization of the PDE is identical to the global problem, with the right hand side populated by the residual from the previous iteration as in equation \eqref{eq:simpleRAS}.

After forming the local active and ghost nodes, there will still be incomplete interpolation and finite difference stencils near the boundary of the subdomain. The nodes needed to complete these stencils are gathered into $\Sigma_j^{BC}$ as active nodes. These are kept separate from $\Sigma_j$ since the extension operator over $\Sigma_j^{BC}$ will be modified in Section \ref{trans} to enforce the transmission conditions.

For Robin transmission conditions, additional ghost nodes need to be layered around the active nodes in $\Sigma_j^{BC}$.  These ghost nodes are also added to $\Sigma_j^{BC}$. The nodes in $\Sigma_j^{BC}$ are shown in Figure \ref{fig:sCons} as diamonds for the active boundary nodes and plus signs for the ghost boundary nodes.

% -------------------------------------------------------------
\subsection{Subsurface representation}
\label{subd:bndGeo}
The active nodes $\Sigma_j$ can be used to partition the surface $\surf$ into a set of overlapping subsurfaces $\surf_j$.
In our discretization of the Robin conditions in Section \ref{trans:robfo}, a closest point representation of the subsurface, $\surf_j$, will be needed.
Conormal vectors, those that are simultaneously orthogonal to the surface normal direction and $\partial\surf_j$, must also be generated to specify the Neumann component of the Robin transmission condition. Note that these constructions are not needed for the simpler Dirichlet boundary conditions presented in Section \ref{trans:dirfo}.

The closest points of the final layer of overlap nodes approximate the subsurface boundary, $\partial\surf_j$. We collect these boundary approximations into a set $\Lambda_j\subset\mathbb{R}^d$.
To obtain a closest point representation for the subsurface, $\surf_j$, we select the closest points for boundary nodes $\Sigma_j^{BC}$ from $\Lambda_j$.  This yields the desired subsurface representation,
\begin{equation}
  CP_{\surf_j}(x_i) =
  \begin{cases}
    CP_\surf(x_i), \quad& x_i\in\Sigma_j\cup\Sigma_j^G, \\
    \underset{y\in \Lambda_j}{\argmin}||x_i-y||_2, \quad& x_i\in\Sigma_j^{BC}.
  \end{cases}
  \label{eq:locCP}
\end{equation}
The closest point computation in the final line is done efficiently by looping over  points in $\Lambda_j$: For a given $y_i\in \Lambda_j$ all boundary nodes within the radius of the interpolation stencil are visited. If a node is unaffiliated, then its closest point is assigned to be $y_i$. Otherwise, if $y_i$ is closer to the node than the previous candidate, then it is updated. This approach to building a closest point function is reminiscent of the method in \cite{CBM:LSE} for constructing closest point functions from triangulated surfaces.

Boundary conormal vectors, \revd{which are simultaneously orthogonal to $\partial\surf_j$ and the surface normal $\hat{n}_i$,} are needed at each of the approximate boundary locations \revd{(see Figure \ref{fig:robfo} for an illustration)}. We approximate the conormal by computing the displacement vector connecting a node $x_i\in\Sigma_j^{BC}$ to its effective boundary location $d_i := x_i - CP_{\surf_j}(x_i)$, and keep only the component that is orthogonal to the surface normal. Upon normalization, a usable approximation to the conormal vector is obtained,
\begin{equation}
  \hat{q}_i = \frac{d_i - \left(d_i\cdot \hat{n}_i\right)\hat{n}_i}{\left||d_i - \left(d_i\cdot \hat{n}_i\right)\hat{n}_i|\right|_2}.
  \label{eq:naConorm}
\end{equation}
Should $d_i$ lie completely in the surface normal direction, $\hat{q}_i$ is set to the zero vector. As is explained in Section \ref{trans:robfo}, \revd{replacing $\hat{q}_i$ in this special case} is robust and corresponds to applying a standard closest point extension from the surface, i.e., it treats the node as if it were an active node in $\Sigma_j$.

% -------------------------------------------------------------
\section{Local problems and transmission conditions}
\label{trans}
The convergence rates of Schwarz type DD solvers and preconditioners depend strongly on the type of transmission conditions posed along the artificial interfaces. In order for the solvers to recover the global discrete solution, the discretization of the transmission conditions must be compatible with the global discretization of the problem \cite{DoleanNataf}. We demonstrate that the forthcoming transmission conditions are compatible with the global discretization through a numerical example in Section \ref{results:transComp}.

This discretization requires special attention for the CPM, as the imposition of boundary conditions occurs through modification of the extension operator prior to composition with the differential operator on the embedding space. As such, generation of transmission operators from splittings of the assembled matrix is not straightforward. For this reason, the subproblems are first considered from the continuous point of view and then discretized, rather than found directly from the algebraic perspective. The local operators will take the form
\begin{equation}
  \mat{A}_j = \left(c+\frac{2d}{h^2}\right)\mat{I} - \left(\frac{2d}{h^2}+\Delta_j^h\right)\begin{bmatrix}\mat{E}_j \\ \mat{T}_j\end{bmatrix},
    \label{eq:localOP}
\end{equation}
where $\mat{E}_j$ is the portion of the global extension operator acting on $\Sigma_j$, and $\mat{T}_j$ is the extension operator, acting over the active nodes in $\Sigma_j^{BC}$, that is modified to enforce the transmission conditions. \revd{As shown in Section \ref{bg:dd}, the linearity of the problem may be used to solve for additive corrections to the solution. This modification also serves to homogenize the transmission conditions.} Accordingly, the final rows of the restricted residual on the right hand side of the local problem, $\mat{A}_je_j=r_j$, are set to zero to match. Throughout this paper, the local problems are solved directly via LU factorization.

% -------------------------------------------------------------
\subsection{Dirichlet transmission conditions}
\label{trans:dirfo}
Homogeneous Dirichlet conditions take the form
\begin{equation}
  u\left(CP_{\surf_j}(x_i)\right) = 0,\quad\forall x_i\in\Sigma_j^{BC},
\end{equation}
and can be enforced to first order accuracy by holding the solution at zero over all nodes in $\Sigma_j^{BC}$, \revd{i.e. by $u(x_i)=0,~\forall x_i\in\Sigma^{BC}_j$,} completely ignoring any underlying geometry. The modified extension is thus the identity over all nodes in $\Sigma_j^{BC}$, and takes the form $\mat{T}_j = \begin{bmatrix} 0 & \mat{I} \end{bmatrix}$ with the zero matrix padding the columns corresponding to interior nodes \revd{$\Sigma_j$}. Since $r_j$ is set to zero in the final entries, corresponding to \revd{the boundary nodes $\Sigma^{BC}_j$}, these transmission conditions yield the same solver as algebraic RAS with a splitting over the same overlapping and disjoint node sets, $\Sigma_j$ and $\widetilde{\Sigma}_j$. This would not be the case for Dirichlet conditions enforced to second order accuracy.

% -------------------------------------------------------------
\subsection{Robin transmission conditions}
\label{trans:robfo}
Optimized restricted additive Schwarz (ORAS) methods enforce Robin transmission conditions. The incorporation of derivative information can dramatically increase the convergence rate of the scheme \cite{DoleanNataf,Gand:OptPar,Cyr:OMSORAS}. Posing these transmission conditions over the same disjoint and overlapping partitions as before leads again to homogeneous boundary conditions when solving for solution corrections.

The discretization of Robin boundary conditions has the side effect of changing the equations that boundary nodes must satisfy, and as such the local operators in \eqref{eq:algRAS} can not be obtained by restrictions of the global operator. This modification of the equations induces a condition on the overlap of the subdomains, requiring the stencil of the Robin transmission operator to lie completely within the overlap region \cite{Cyr:OMSORAS}. This, paired with the relative lack of importance of the order of accuracy of the imposed transmission conditions, motivates the construction of a Robin boundary operator for the CPM using a minimal stencil size. We propose a first order accurate discretization here that combines a forward difference for the Neumann part with a point evaluation for the Dirichlet part.

\begin{figure}[htbp]
  \centering
  \includegraphics[width=0.75\linewidth]{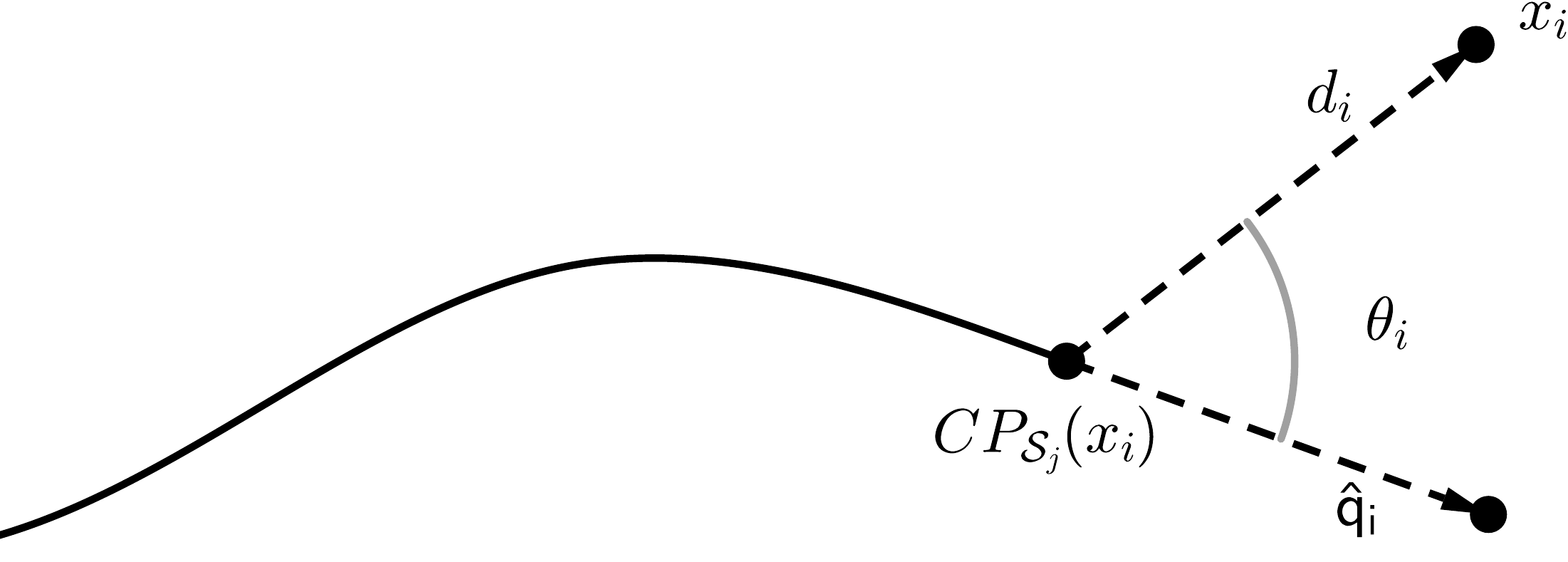}
  \caption{Here the relevant components for the first-order accurate Robin condition are shown. The boundary node \revd{$x_i\in\Sigma^{BC}_j$} is associated with its closest point on the surface $CP_{\surf_j}(x_i)$, from which the displacement vector $d_i$ and conormal vector $\hat{q}_i$ are constructed.}
  \label{fig:robfo}
\end{figure}

The Robin condition
\begin{equation}
  \restrict{\pdiff{u}{\hat{q}_i}}{CP_{\surf_j}(x_i)} + \alpha u\left(CP_{\surf_j}(x_i)\right) = 0,
  \label{eq:contRob}
\end{equation}
will be enforced at each boundary location \revd{$CP_{\surf_j}(x_i)\in\Lambda_j$}; see Figure \ref{fig:robfo} for an illustration. Notice that the partial derivative in the boundary conormal direction, $\hat{q}_i$, can be written in terms of the node's displacement from its effective boundary position through a simple change of variables. This yields
\begin{equation}
  \restrict{\pdiff{u}{d_i}}{CP_{\surf_j}(x_i)} = \cos\theta_i\restrict{\pdiff{u}{\hat{q}_i}}{CP_{\surf_j}(x_i)} + \sin\theta_i\restrict{\pdiff{u}{\hat{n}_i}}{CP_{\surf_j}(x_i)},
  \label{eq:pdiffChgVar}
\end{equation}
where $\theta_i$ is the angle between $\hat{q}_i$ and $d_i$. The partial derivative in the surface normal direction vanishes since the solution is constant along the surface normals. Assuming for a moment that $\hat{q}_i$ and $d_i$ are not perpendicular, we may combine  equations \eqref{eq:contRob} and \eqref{eq:pdiffChgVar} to find
\begin{equation}
  \frac{1}{\cos\theta_i}\restrict{\pdiff{u}{d_i}}{CP_{\surf_j}(x_i)} + \alpha u\left(CP_{\surf_j}(x_i)\right) = 0,
  \label{eq:contRobChgVar}
\end{equation}
which is easier to discretize. Replacing the partial derivative with a forward difference and discretizing the Dirichlet term as in Section~\ref{trans:dirfo} yields
\begin{equation}
  \frac{u(x_i)-u\left(CP_{\surf_j}(x_i)\right)}{d_i\cdot\hat{q}_i} + \alpha u(x_i) = 0,
  \label{eq:discRob}
\end{equation}
where the angle $\theta_i$ has been absorbed into the dot product $d_i\cdot\hat{q}_i$.
Finally, isolating $u(x_i)$ in this expression gives the relationship that the extension operator must satisfy over the boundary nodes, \revd{$\Sigma^{BC}_j$}, as
\begin{equation}
  u(x_i) = \frac{u\left(CP_{\surf_j}(x_i)\right)}{1 + \alpha d_i\cdot\hat{q}_i}.
  \label{eq:robfo}
\end{equation}
As $d_i\cdot\hat{q}_i\rightarrow 0$, Equation \eqref{eq:robfo} reduces to $u(x_i)=u\left(CP_{\surf_j}(x_i)\right)$ which is identical to the extension step applied to the active nodes $\Sigma_j$. This makes sense, as this case only arises when a boundary node is adjacent to an active node. Thus the transmission condition \eqref{eq:robfo} remains robust against the apparent possibility of division by zero in \eqref{eq:discRob}.

% ---------------------------------------------------------------
\subsubsection{Order of accuracy}
\label{trans:robfo:acc}
We test this discretization \revd{of Robin conditions} on a global problem with a known exact solution. Consider $\surf$ to be a 2 radian arc of the unit circle. The boundary value problem
\begin{equation}
  \begin{cases}
    \left(1-\Delta_\surf\right)u = 1-\theta^2\quad&{\rm on}~\surf,\\
    u(0) = 0,\\
    \restrict{\pdiff{u}{\theta}}{\theta=2} + u(2) = 0,
  \end{cases}
  \label{eq:robConv}
\end{equation}
has the solution \revd{$u(\theta) = \cosh(\theta) + (9e^{-2}-1)\sinh(\theta) - \theta^2 - 1$}. The Laplace-Beltrami operator is discretized via the CPM described in Section \ref{bg:cpm} with bi-cubic interpolation for the extension operator and a second-order centered difference Laplacian on the embedding space. To isolate the effect of the Robin condition, the Dirichlet condition at $\theta=0$ is enforced to second order accuracy by the mirror point discretization described in \cite{CBM:Eig}. Finally, the Robin condition at $\theta=2$ is enforced by the discretization given in equation \eqref{eq:robfo}. Notice that placing the Robin condition at $\theta=2$ yields a conormal direction that is not aligned with the grid on the embedding space. The error in the solution is measured in the infinity norm over a sequence of grid resolutions, as seen in Table \ref{tab:robConv}. We observe that the global error is first-order accurate.

\begin{table}[hbtp]
  \centering \small
  \begin{tabular}{c|ccccc}
  %  \toprule
    $h$                   &  $1/64$ & $1/128$ & $1/256$ & $1/512$ & $1/1024$\\ \midrule
    $||u-u_{ex}||_\infty$ &  $5.15\times 10^{-2}$ & $2.56\times 10^{-2}$ & $1.27\times 10^{-2}$ & $6.36\times 10^{-3}$ & $3.18\times 10^{-3}$\\ %\midrule
%    $h$                   & $1/256$ & $1/512$ & $1/1024$ \\
 %   $||u-u_{ex}||_\infty$ & $1.27\times 10^{-2}$ & $6.36\times 10^{-3}$ & $3.18\times 10^{-3}$ \\ \bottomrule
  \end{tabular}
  \caption{The global $\infty-$norm error for a CPM discretization of an elliptic problem on a circular arc \eqref{eq:robConv}. The problem was solved on 5 grids of increasing refinement using the proposed Robin boundary condition discretization~\eqref{eq:robfo}. We observe first order accuracy.}
  \label{tab:robConv}
\end{table}

% ---------------------------------------------------------------
\subsubsection{Modified scheme for cross points}
\label{trans:robfo:crosspt}
Points where more than two subdomains meet are anticipated when a surface embedded in $\mathbb{R}^d$, $d\geq 3$, is partitioned. Non-overlapping optimized Schwarz methods are destabilized by these points when otherwise optimal values of the Robin parameter are used \cite{Gand:XPT, Loisel:2Lag}. This behavior typically does not appear in overlapping optimized Schwarz (ORAS) methods, used as solvers or as preconditioners, since each unknown is only updated by a single subproblem. However, we do see this behavior for the CPM, and attribute it to the presence of the extension operator enlarging the stencil size.

Stability and convergence can be recovered by increasing $\alpha$, thereby weighting the Dirichlet component of the Robin transmission condition more heavily. However, raising the value of $\alpha$ globally is non-optimal; fewer iterations are required if we use the larger value only in the vicinity of the cross points. This approach weighs the Dirichlet term more heavily in regions where it is needed for stability, while selecting weights elsewhere for rapid convergence. We let $\alpha^\times$ denote the larger weight used near the cross points, and leave $\alpha$ to denote the Robin weight elsewhere. For Poisson and positive Helmholtz problems on the plane it is known that the optimal values of these weights scale as $\alpha\sim\mathcal{O}\left(h^{-1/2}\right)$ away from the cross points and $\alpha^\times\sim\mathcal{O}\left(h^{-1}\right)$ near the cross points \cite{Gand:XPT}.

Identifying regions for applying each weight proceeds in a fashion similar to the construction of the boundary locations in Section \ref{subd:bndGeo}. Cross point nodes are identified as those in each disjoint subdomain $\widetilde{\Sigma}_j$ that have neighbors belonging to more than one other subdomain. A sphere of radius $2N_Oh$ is centered on each of these cross point nodes, and all boundary nodes in $\Sigma_j^{BC}$ lying within this sphere are marked to use a Robin weight of $\alpha^\times$ in place of $\alpha$.

% -------------------------------------------------------------
\section{Results}
\label{results}
Our development of (O)RAS methods for the iterative solution of the stabilized CPM has seen the introduction of the standard DD parameters, \revd{the number of subdomains} $N_S$, \revd{the overlap width} $N_O$, and \revd{Robin weight} $\alpha$. \revd{A} subdomain boundary alignment method and two transmission conditions specific to this problem have also been introduced. In the following sub-sections we sweep over these parameters and evaluate the effect of these problem-dependent choices. The software written to solve these problems is available through a public repository \cite{May:CPMCode}. These results were obtained specifically with commit \revd{\textit{fc3612f9}}. Later commits will be focused on user friendliness and general applicability of the code-base.

\begin{figure}[hbtp]
  \centering
  \includegraphics[width=0.75\linewidth]{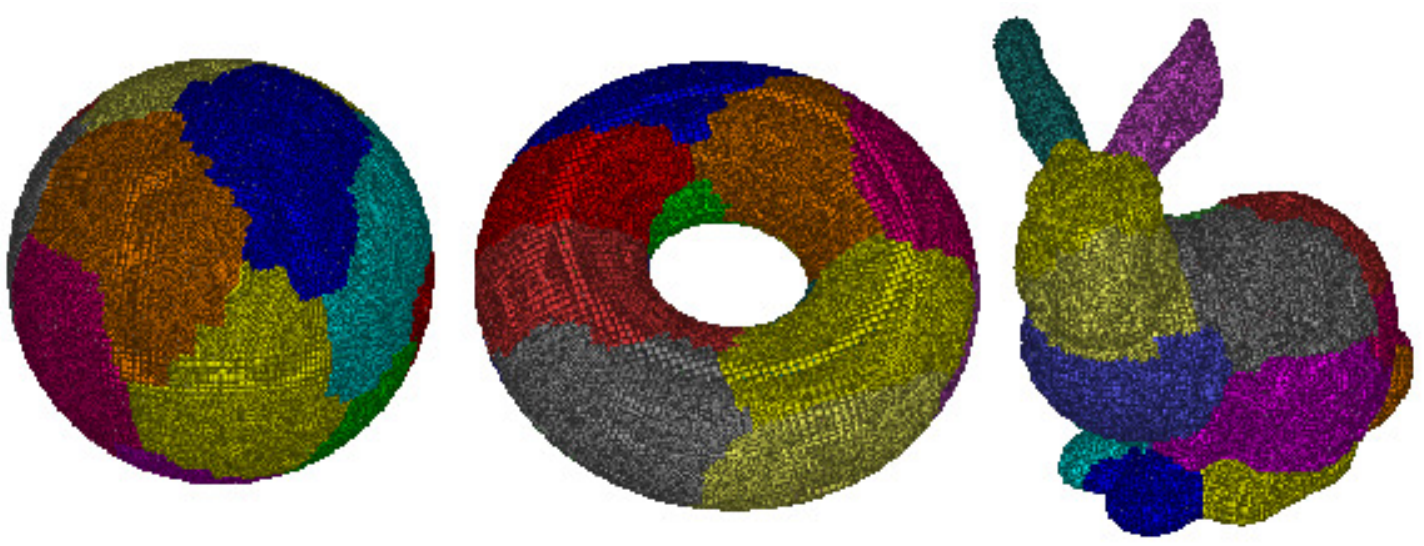}
  \caption{The three surfaces used to evaluate our solvers with colors showing individual subdomains obtained from METIS. All surfaces use $h=1/30$ and have 16 subdomains. The active nodes are projected onto the surfaces for the sake of visualization.}
  \label{fig:exSurfs}
\end{figure}

The methods are tested over three surfaces of increasing complexity, as seen in Figure \ref{fig:exSurfs}. The unit sphere is simple, yet extremely important, and makes an ideal first test. The torus, with major radius $R=2/3$ and minor radius $r=1/3$, introduces a hole and regions of negative curvature. Finally, the Stanford Bunny \cite{bunny} is a triangulated surface with regions of high curvature and holes. We use the original Stanford Bunny triangulation with 5 holes. Our only modification is to scale the bunny to be 2 units tall. In all cases, the extension operator is constructed with tri-quadratic tensor product interpolation, as introduced in Section \ref{bg:cpm}. In all tests, convergence was declared when the $2-$norm of the residual vector was reduced by a factor of $10^6$.

% -------------------------------------------------------------
\subsection{Consistency of the transmission conditions}
\label{results:transComp}
Before evaluating the performance of the proposed solvers and the effects of the various DD parameters, the transmission conditions defined in Section \ref{trans} are validated for compatibility with the global discretization. This is done numerically by comparing the DD solution to the direct solution of the global system. We consider the unit sphere with $h=1/25$ and $h=1/50$, divided into $N_S=16$ subdomains with an overlap of $N_O=4$. The Robin transmission conditions use the weights $\alpha=4$ and $\alpha^\times=40$. The reference solution of the global system is obtained with the sparse direct solver MUMPS \cite{MUMPS01,MUMPS02} for each resolution. The right hand side is chosen such that the exact continuous solution is $\sin^2(\phi)e^{\cos\theta}$, to verify that the error between the DD solution and global discrete solution is below the truncation error present in the discretization.

\begin{figure}[hbtp]
  \centering
  \includegraphics[width=0.85\linewidth]{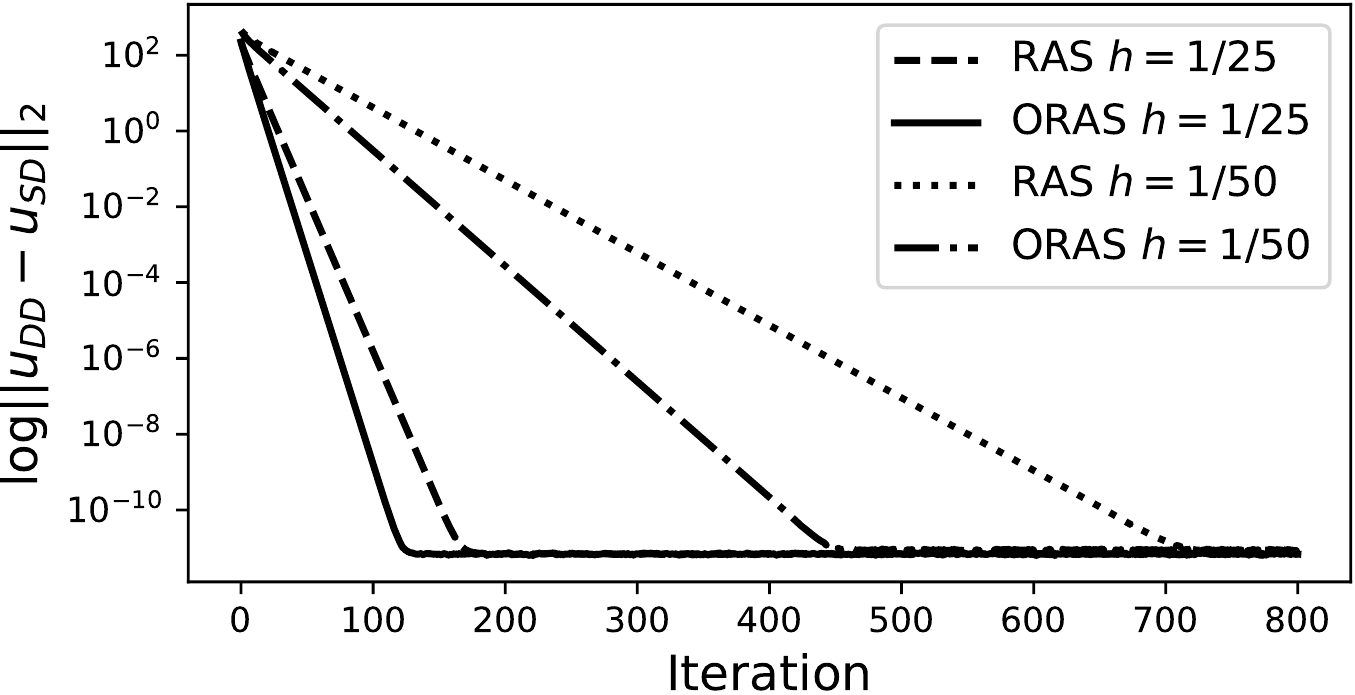}
  \caption{The norm of the difference between the DD solution, $U_{DD}$, and the MUMPS reference solution, $U_{SD}$, is shown here as a function of the DD iteration number for RAS and ORAS solvers on the unit sphere. Grid resolutions of $h=1/25$ and $h=1/50$ were tested. The error compared to the exact continuous solution is $4.486\cdot 10^{-3}$ for $h=1/25$, and $1.264\cdot 10^{-3}$ for $h=1/50$, measured in the $2-$norm.}
  \label{fig:convHist}
\end{figure}

Figure \ref{fig:convHist} shows the $2-$norm of the difference between the reconstructed DD solution and the MUMPS reference solution for each combination of resolution and (O)RAS solver. Notably, the error in the DD solution plateaus at a significantly smaller value than the truncation error present in the discretization itself. Indeed, the error plateau only depends on the condition number of the global matrix, and thus on the reliability of the reference MUMPS solution. No dependence on the choice of solver, right hand side, or DD parameters has been observed.

% -------------------------------------------------------------
\subsection{Effect of the Robin weight}
\label{results:alpha}
Next, we investigate how the optimal Robin parameter, $\alpha$, depends on the grid spacing for the three surfaces by splitting each into two subdomains with a fixed overlap of $N_O=4$ and varying $\alpha$ and $h$. There are no cross points when the surface is split into two subdomains.

The unit sphere produces meshes with $N_A=163542$, $N_A=654630$ and $N_A=2618046$ for grid spacings of $h=1/50$, $h=1/100$, and $h=1/200$, respectively. The torus produces meshes with $N_A=114692$, $N_A=454300$ and $N_A=1818052$ for the same grid spacings. Finally, the Stanford Bunny produces meshes with $N_A=98714$, $N_A=396232$ and $N_A=1585218$, respectively.

\begin{table}[htb]
  \centering
  \begin{tabular}{c|cccccccccc}
    \toprule
    \diagbox{$h$}{$\alpha$} & $1/4$ & $1/2$ & $3/4$ & 1  & $3/2$ & 2  & 4  & 8  & 16 & $\infty$ \\ \midrule
    & \multicolumn{10}{c}{Solver} \\
    $1/50$   & 19    & 14 & 11 & 11 & 15 & 20 & 33 & 49 & 67  & 83   \\
    $1/100$  & 25    & 17 & 14 & 13 & 16 & 21 & 37 & 61 & 93  & 159  \\
    $1/200$  & 32    & 22 & 18 & 18 & 18 & 23 & 43 & 75 & 122 & 303  \\
    & \multicolumn{10}{c}{Preconditioner} \\
    $1/50$   & 11    & 9  & 9  & 8  & 8  & 9  & 11 & 13 & 15  & 17 \\
    $1/100$  & 12    & 11 & 10 & 9  & 9  & 9  & 12 & 15 & 18  & 20 \\
    $1/200$  & 14    & 12 & 12 & 11 & 10 & 11 & 13 & 17 & 21  & 31 \\ \bottomrule
  \end{tabular}
  \caption{Iterations to convergence for the unit sphere as the Robin \revd{weight} $\alpha$ and grid spacing $h$ were varied.  The sphere was divided into two subdomains with an overlap of $N_O=4$. The final column gives the corresponding RAS iteration count.}
  \label{tab:sphTwoSub}
\end{table}

Table \ref{tab:sphTwoSub} shows the number of iterations to convergence for the sphere with grid resolutions of $h=1/50,1/100,1/200$ as $\alpha$ is varied. Tables \ref{tab:torTwoSub} and \ref{tab:triTwoSub} show the same information for the torus and the Stanford Bunny, respectively. The ORAS solvers on the sphere and torus show flat iteration counts around the optimal $\alpha$ value. Though $\alpha=1$ appears optimal throughout all resolutions for these two surfaces, we note that the iteration counts for $\alpha=1/2$ are growing  and the region of flat iteration counts is shifting to larger values of $\alpha$ as the resolution is refined. The ORAS solver on the Stanford Bunny shows a clear dependence of the optimal value for $\alpha$ on the resolution, varying from $\alpha=1$ when $h=1/50$ to $\alpha=2$ when $h=1/200$.
For all three surfaces, the ORAS preconditioner shows the same trend, with the optimal $\alpha$ value increasing with grid refinement, but the effect is weaker.

In each of these cases, the iteration counts of the ORAS methods are seen to approach the iteration count of the corresponding RAS method as $\alpha$ is increased beyond the optimal value. This agrees with our expectation, since the weight of the Dirichlet term increases as $\alpha$ increases.

\begin{table}[htb]
  \centering
  \begin{tabular}{c|cccccccccc}
    \toprule
    \diagbox{$h$}{$\alpha$} & $1/4$  & $1/2$ & $3/4$ & 1 & $3/2$  & 2  & 4  & 8   & 16 & $\infty$ \\ \midrule
    & \multicolumn{10}{c}{Solver} \\
    $1/50$   & 20     & 13 & 11 & 10 & 14 & 18 & 30 & 44  & 59  & 74  \\
    $1/100$  & 24     & 17 & 13 & 12 & 16 & 20 & 36 & 58  & 86  & 142 \\
    $1/200$  & 29     & 20 & 17 & 15 & 15 & 20 & 37 & 66  & 107 & 272 \\
    & \multicolumn{10}{c}{Preconditioner} \\
    $1/50$   & 10     & 9  & 8  & 7  & 7  & 8  & 9  & 10  & 12  & 13 \\
    $1/100$  & 11     & 10 & 9  & 8  & 8  & 8  & 10 & 12  & 14  & 17 \\
    $1/200$  & 12     & 11 & 10 & 9  & 9  & 9  & 10 & 12  & 15  & 22 \\ \bottomrule
  \end{tabular}
  \caption{Iterations to convergence for the torus as the Robin \revd{weight} $\alpha$ and grid spacing $h$ are varied. The torus was divided into two subdomains with an overlap of $N_O=4$. The final column gives the corresponding RAS iteration count.}
  \label{tab:torTwoSub}
\end{table}

\begin{table}[htb]
  \centering
  \begin{tabular}{c|cccccccccc}
    \toprule
    \diagbox{$h$}{$\alpha$} & $1/4$  & $1/2$ & $3/4$ & 1 & $3/2$ & 2  & 4  & 8   & 16 & $\infty$ \\ \midrule
    & \multicolumn{10}{c}{Solver} \\
    $1/50$   & 18     & 14 & 11 & 11 & 15 & 19 & 30 & 42  & 55  & 65  \\
    $1/100$  & 69     & 30 & 21 & 17 & 16 & 20 & 35 & 57  & 83  & 133 \\
    $1/200$  & 46     & 33 & 27 & 24 & 20 & 18 & 32 & 55  & 89  & 215 \\
    & \multicolumn{10}{c}{Preconditioner} \\
    $1/50$   & 11     & 9  & 8  & 8  & 9  & 9  & 10 & 12  & 13  & 14 \\
    $1/100$  & 15     & 12 & 11 & 10 & 10 & 10 & 10 & 14  & 16  & 20 \\
    $1/200$  & 16     & 15 & 14 & 13 & 12 & 11 & 12 & 15  & 18  & 28 \\ \bottomrule
  \end{tabular}
  \caption{Iterations to convergence for the Stanford Bunny as the Robin \revd{weight} $\alpha$ and grid spacing $h$ were varied. The bunny was divided into two subdomains with an overlap of $N_O=4$. The final column gives the corresponding RAS iteration count. The optimal value of $\alpha$ grows as the grid resolution is refined.}
  \label{tab:triTwoSub}
\end{table}

% -------------------------------------------------------------
\subsection{Necessity of the cross point modification}
\label{results:alphacross}
In this sub-section, we will examine two resolutions of the grid, $h = 1/100$ and $h = 1/200$, with the same number of global unknowns as the corresponding problems in the previous section. In the preceding section, only two subdomains were used, and thus there were no cross points present in the decomposition. Here we decompose the problems into $N_S=32$ subdomains, which introduces cross points into the splitting. The Robin weights $\alpha$ and $\alpha^\times$ are varied and the required iterations to convergence for the ORAS solvers and preconditioners gathered. With cross points present, the values of $\alpha$ permitting convergence are shifted upwards, and in the following $\alpha=2$ is the minimum tested.

\begin{figure}[hbtp]
  \centering
  \includegraphics[width=0.75\linewidth]{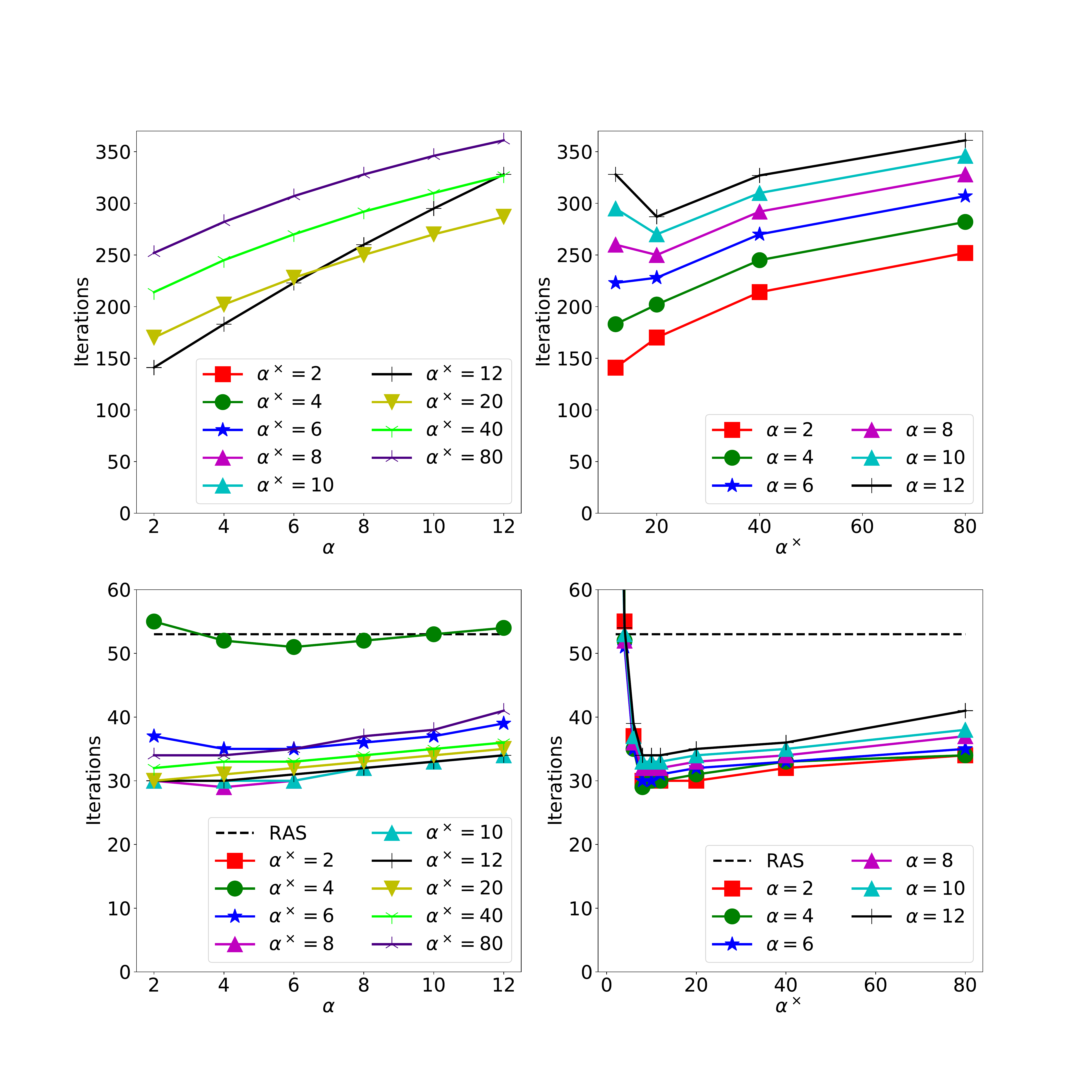}
  \caption{Dependence of iterations to convergence on $\alpha$ and $\alpha^\times$ for the sphere with $h=1/100$. The upper panels show the iterations to convergence when ORAS is run as a solver, while the lower panels use ORAS as a preconditioner. The left panels sweep over different values of $\alpha$, holding $\alpha^\times$ fixed. The right panels show the same results as a function of $\alpha^\times$ for fixed values of $\alpha$. As shown in the upper left panel, $\alpha^\times<12$ gives solvers that never converge regardless of the choice of $\alpha$. For comparison, RAS requires $500$ iterations as a solver and $53$ iterations as a preconditioner.}
  \label{fig:sphereAlphaCross}
\end{figure}

\begin{figure}[hbtp]
  \centering
  \includegraphics[width=0.75\linewidth]{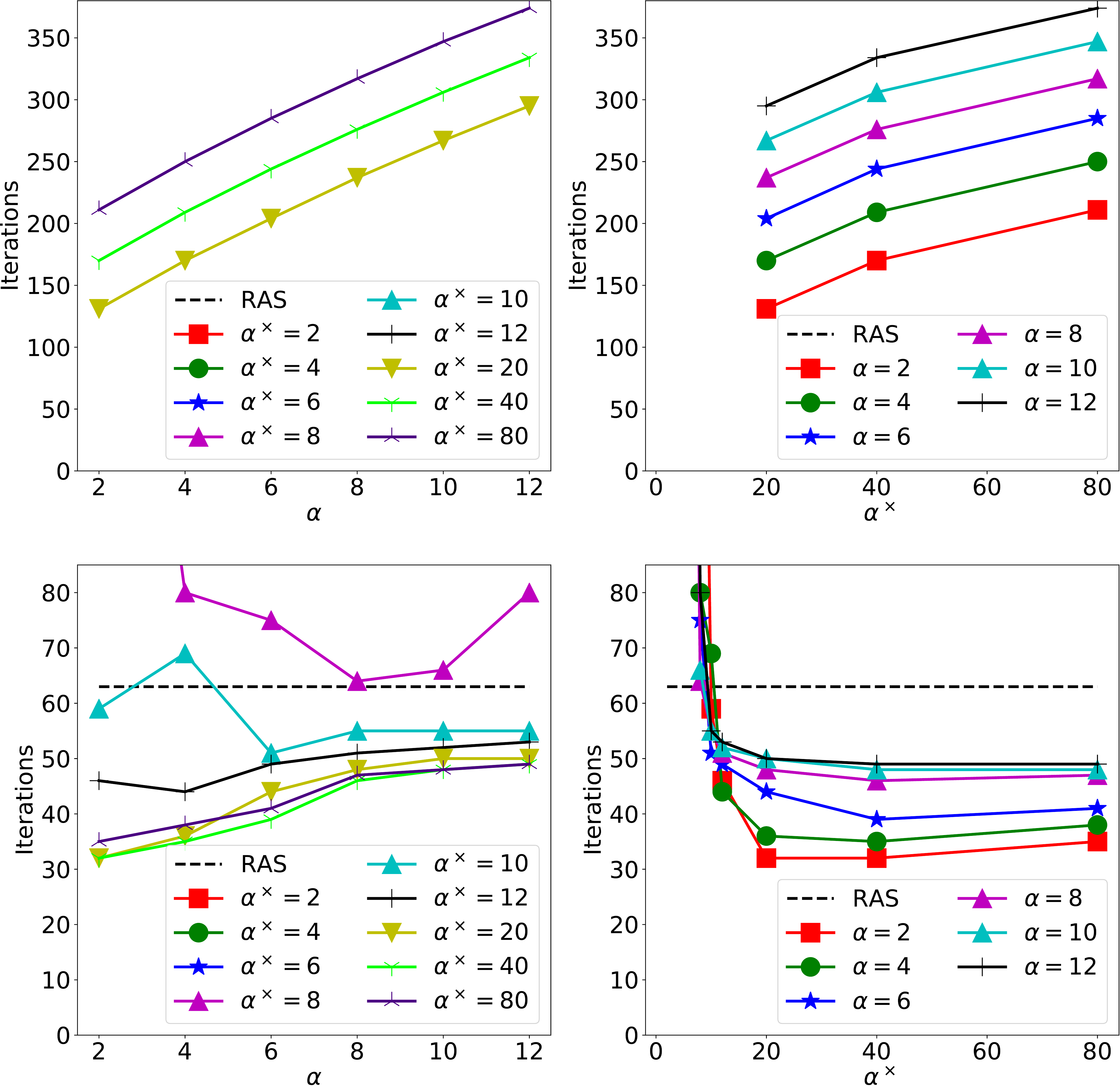}
  \caption{Dependence of iterations to convergence on $\alpha$ and $\alpha^\times$ for the sphere with $h=1/200$. The upper panels show the iterations to convergence when ORAS is run as a solver, while the lower panels use ORAS as a preconditioner. The left panels sweep over different values of $\alpha$, holding $\alpha^\times$ fixed. The right panels show the same results as a function of $\alpha^\times$ for fixed values of $\alpha$. As shown in the upper left panel, $\alpha^\times<20$ gives solvers that never converge regardless of the choice of $\alpha$. For comparison, RAS requires $916$ iterations as a solver and $63$ iterations as a preconditioner.}
  \label{fig:sphereAlphaCrossHR}
\end{figure}

We consider the sphere first. In Figures \ref{fig:sphereAlphaCross} and \ref{fig:sphereAlphaCrossHR} we see that the required iterations of the solver and preconditioner depend more strongly on the value of $\alpha$ than $\alpha^\times$, provided $\alpha^\times$ is large enough to yield a convergent method. In particular, when $h=1/100$, $\alpha=2$, and $\alpha^\times$ is set to $12$, $20$, and $40$, the solver (preconditioner) requires $141(30)$, $170(30)$, and $214(32)$ iterations, respectively. Similarly for $h=1/200$, setting $\alpha=2$, and $\alpha^\times$ to $20$, $40$, and $80$, gives solvers (preconditioners) that require $131(32)$, $170(32)$, and $211(35)$ iterations to converge.

\begin{figure}[hbtp]
  \centering
  \includegraphics[width=0.75\linewidth]{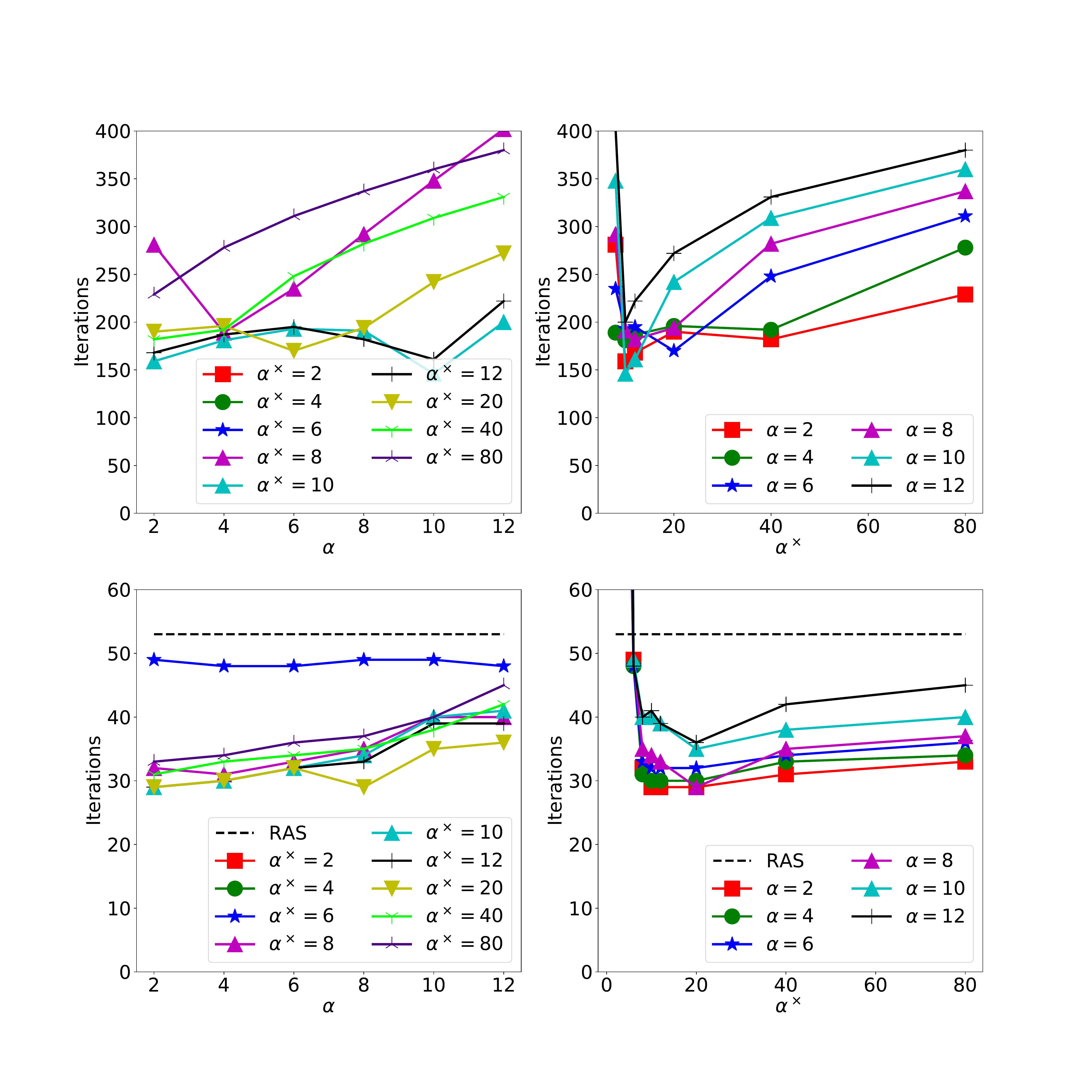}
  \caption{Dependence of iterations to convergence on $\alpha$ and $\alpha^\times$ for the torus with $h=1/100$. The upper panels show the iterations to convergence when ORAS is run as a solver, while the lower panels use ORAS as a preconditioner. The left panels sweep over different values of $\alpha$, holding $\alpha^\times$ fixed. The right panels show the same results as a function of $\alpha^\times$ for fixed values of $\alpha$. As shown in the upper left panel, $\alpha^\times<8$ gives solvers that never converge regardless of the choice of $\alpha$. For comparison, RAS requires $597$ iterations as a solver and $53$ iterations as a preconditioner.}
  \label{fig:torusAlphaCross}
\end{figure}

\begin{figure}[htbp]
  \centering
  \includegraphics[width=0.75\linewidth]{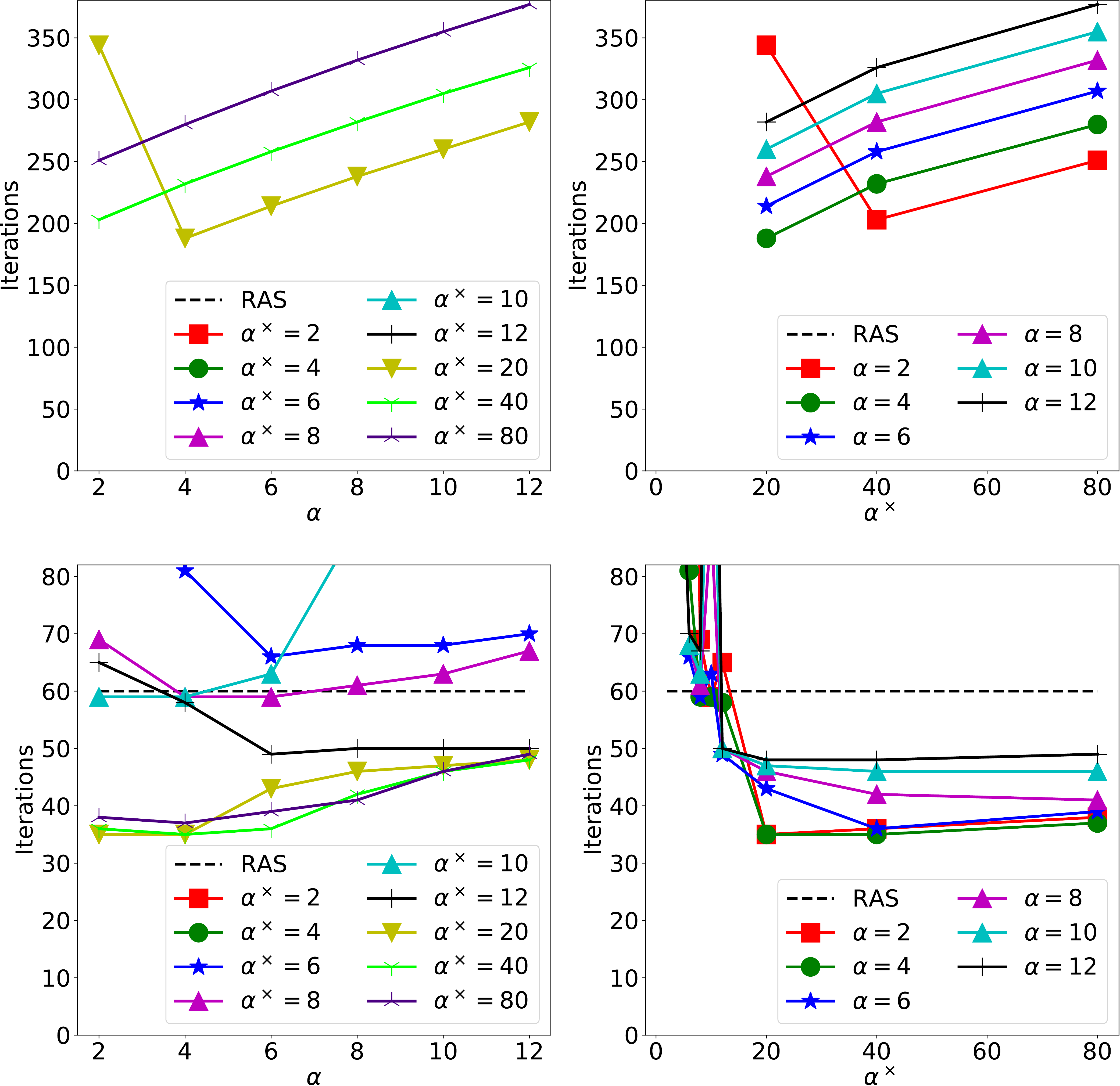}
  \caption{Dependence of iterations to convergence on $\alpha$ and $\alpha^\times$ for the torus with $h=1/200$. The upper panels show the iterations to convergence when ORAS is run as a solver, while the lower panels use ORAS as a preconditioner. The left panels sweep over different values of $\alpha$, holding $\alpha^\times$ fixed. The right panels show the same results as a function of $\alpha^\times$ for fixed values of $\alpha$. As shown in the upper left panel, $\alpha^\times<20$ gives solvers that never converge regardless of the choice of $\alpha$. For comparison, RAS requires $1148$ iterations as a solver and $60$ iterations as a preconditioner.}
  \label{fig:torusAlphaCrossHR}
\end{figure}

Figures \ref{fig:torusAlphaCross} and \ref{fig:torusAlphaCrossHR} show the same sweeps over $\alpha$ and $\alpha^\times$, but for the torus. For $h = 1/100$ there is an optimum parameter combination for the ORAS solver with $\alpha=10$ and $\alpha^\times=10$, requiring only $146$ iterations for the solver to converge. Although this implies that the cross point modification is unnecessary, it is the only time this was observed and seems to be coincidental with the splitting produced by METIS for this surface at this resolution and subdomain count. At $h = 1/200$ the optimal solver uses $\alpha=4$ and $\alpha^\times=20$. This decrease in the optimal value for $\alpha$ opposes the expected behavior, although the large number of subdomains and complicated geometry make comparisons to the flat case difficult. At this resolution $\alpha^\times=20$ is optimal for almost all values of $\alpha$. Consistent with the results for the sphere, the dependence of the preconditioners on $\alpha^\times$ is weak so long as the value is large enough to yield a convergent solver.

\begin{figure}[htbp]
  \centering
  \includegraphics[width=0.75\linewidth]{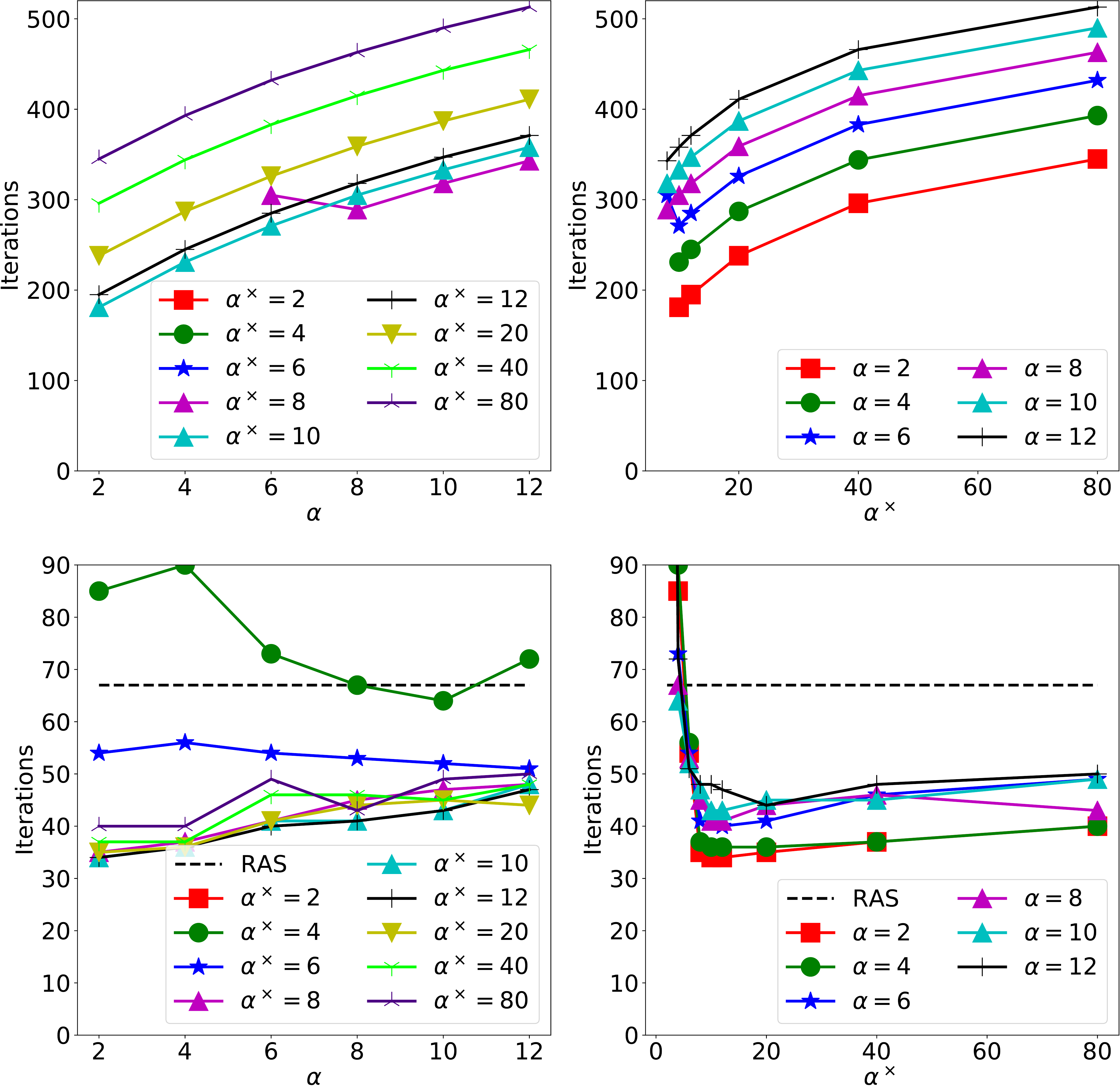}
  \caption{Dependence of iterations to convergence on $\alpha$ and $\alpha^\times$ for the Stanford Bunny with $h=1/100$. The upper panels show the iterations to convergence when ORAS is run as a solver, while the lower panels use ORAS as a preconditioner. The left panels sweep over different values of $\alpha$, holding $\alpha^\times$ fixed. The right panels show the same results as a function of $\alpha^\times$ for fixed values of $\alpha$. As shown in the upper left panel, $\alpha^\times<8$ gives solvers that never converge regardless of the choice of $\alpha$. For comparison, RAS requires $717$ iterations as a solver and $67$ iterations as a preconditioner.}
  \label{fig:triangAlphaCross}
\end{figure}

\begin{figure}[htbp]
  \centering
  \includegraphics[width=0.75\linewidth]{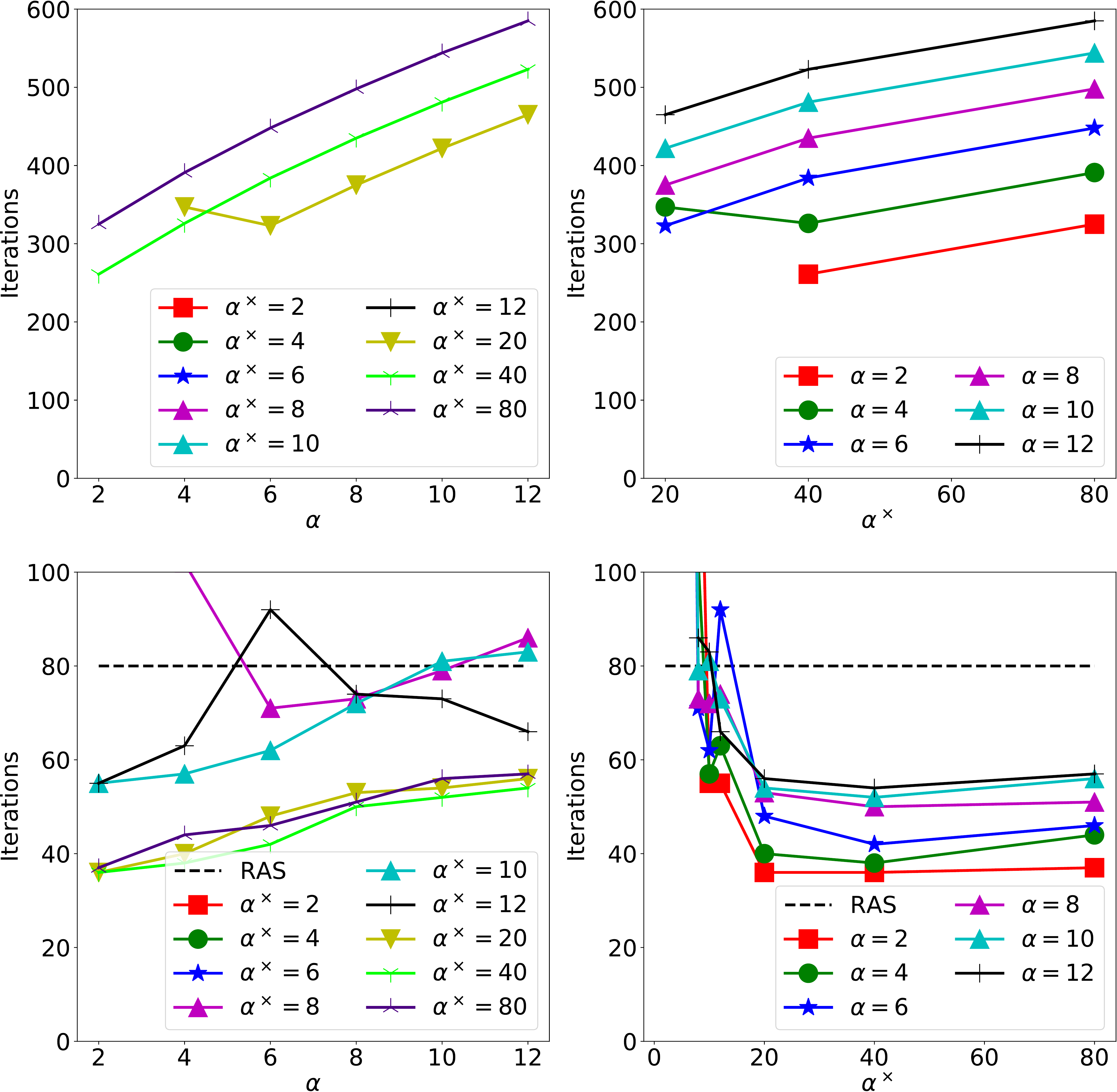}
  \caption{Dependence of iterations to convergence on $\alpha$ and $\alpha^\times$ for the Stanford Bunny with $h=1/200$. The upper panels show the iterations to convergence when ORAS is run as a solver, while the lower panels use ORAS as a preconditioner. The left panels sweep over different values of $\alpha$, holding $\alpha^\times$ fixed. The right panels show the same results as a function of $\alpha^\times$ for fixed values of $\alpha$. For comparison, RAS requires $1301$ iterations as a solver and $80$ iterations as a preconditioner.}
 \label{fig:triangAlphaCrossHR}
\end{figure}

Finally, Figures \ref{fig:triangAlphaCross} and \ref{fig:triangAlphaCrossHR} show the effects of $\alpha$ and $\alpha^\times$ on the iterations to convergence for the ORAS solvers and preconditioners on the Stanford Bunny. Overall the results are quite similar to the sphere and torus problems.

\revd{Notably, the use of of larger values of $\alpha^\times$ allows smaller values of $\alpha$ to be used successfully. For instance, Figure \ref{fig:triangAlphaCrossHR} shows that choosing $\alpha^\times=40$ allows $\alpha=2$ to be used, which, in combination provides the smallest number of iterations to convergence. Conversely, taking $\alpha^\times = 12$ or smaller leads slow or failed convergence (for example, the ORAS solvers did not converge at all for these values). We further observe that the use of the cross point modification with the smaller allowable values of $\alpha$ leads to ORAS solvers and preconditioners that enjoy greatly reduced iteration counts relative to their RAS counterparts.} The preconditioners again depend very weakly on the value of $\alpha^\times$ for values that are sufficiently large.

% -------------------------------------------------------------
\subsection{Effect of the overlap width}
\label{results:nover}
The effect of the overlap width was investigated with $N_S=32$ and $h=1/200$ fixed for all surfaces. The ORAS methods use Robin parameters of $\alpha=2$ and $\alpha^\times=40$ for the sphere, $\alpha=4$ and $\alpha^\times=20$ for the torus, and $\alpha=2$ and $\alpha^\times=40$ for the Stanford Bunny.

\begin{table}[htbp]
  \centering
  \begin{tabular}{@{}r|cc|cc|cc@{}}
    \toprule
    $N_O$            & \multicolumn{2}{c|}{4} & \multicolumn{2}{c|}{8} & \multicolumn{2}{c}{12}  \\
    Method           & RAS        & ORAS       & RAS        & ORAS       & RAS        & ORAS       \\ \midrule
    \multicolumn{7}{c}{Sphere, $h=1/200$, $N_A=2618046$} \\
    Solver         &     916    &    170      &    598     &     182     &     449    &    207      \\
    Preconditioner &  63        & 36          & 55         &   32        &  51        & 31          \\ \midrule
    \multicolumn{7}{c}{Torus, $h=1/200$, $N_A=1818052$} \\
    Solver         &    738     &     188     &    536     &     217     &     412    &    240      \\
    Preconditioner &  60        &  35         & 49         &  33         &  49        & 32          \\ \midrule
    \multicolumn{7}{c}{Stanford Bunny, $h=1/200$, $N_A=1585218$} \\
    Solver         & 1301       & 347         & 820        & 307         & 609        & 317         \\
    Preconditioner & 80         & 40          & 64         & 36          & 58         & 35          \\ \bottomrule
  \end{tabular}
  \caption{The iterations to convergence as the overlap width $N_O$ \revd{was} varied \revd{with cross point modification}. The RAS solvers and preconditioners consistently require less iterations to converge as $N_O$ grows, as expected. The ORAS solvers defy expectation and show growing iterations to convergence as $N_O$ grows. This is due to cross point modification affecting more nodes at larger overlap widths. This is confirmed by Table \ref{tab:overAllNC} where the same sweep is done without cross point modification.}
  \label{tab:overAll}
\end{table}

\begin{table}[htbp]
  \centering
  \begin{tabular}{r|ccc|ccc|ccc}
    \toprule
    & \multicolumn{3}{c|}{Sphere} & \multicolumn{3}{c|}{Torus} & \multicolumn{3}{c}{Stanford Bunny} \\
    $N_O$          & 4   & 8   & 12  & 4   & 8   & 12  & 4   & 8   & 12   \\ \midrule
    Solver         & 393 & 319 & 299 & 358 & 312 & 286 & 604 & 503 & 426  \\
    Preconditioner & 53  & 49  & 38  & 51  & 47  & 43  & 64  & 53  & 44   \\ \bottomrule
  \end{tabular}
  \caption{The iterations to convergence as the overlap width $N_O$ was varied without cross point modification. As $N_O$ is increased, the ORAS solvers and preconditioners show consistently decreasing iterations to convergence. Although this matches the expectation for iteration counts against $N_O$, Table \ref{tab:overAll} shows that usage of the cross point modification yields faster convergence regardless of overlap width.}
  \label{tab:overAllNC}
\end{table}

As seen in Table \ref{tab:overAll}, the RAS solvers and preconditioners consistently require fewer iterations as the overlap is widened. The ORAS solvers, contrary to expectation, show increasing iteration counts with widened overlaps while the preconditioners show little variation. Referring to Section \ref{trans:robfo:crosspt}, the cross point modification scheme is dependent on the overlap width. Thus as $N_O$ is increased, the number of affected nodes grows and $\alpha^\times$ is used more often. This shifts the method away from optimal behavior and yields larger iteration counts. To affirm this, the same surfaces and splittings were run without the cross point modification. The results of this, shown in Table \ref{tab:overAllNC}, show decreasing iteration counts as the overlap widens. However, the larger value of $\alpha$, now needed everywhere, slows all of the solvers and preconditioners. For ORAS methods, it is wise to use a small overlap with appropriately chosen parameters $\alpha$ and $\alpha^\times$.

% -------------------------------------------------------------
\subsection{Effect of the subdomain count}
\label{results:nsub}
As the number of subdomains is increased, the subproblems become cheaper to assemble and factor, and the method exhibits greater parallelism. On the other hand, the local problems see less of the global problem, and naturally require more iterations for convergence.

We now examine this in greater detail on our three shapes using a fixed overlap of $N_O=4$, and subdomain counts $N_S=32,64,96$, with grid spacings of $h=1/100$, and $1/200$. With these large subdomain counts, there are many cross points scattered through the domain. The cross point modified scheme, introduced in Section \ref{trans:robfo:crosspt} and examined in Section \ref{results:alphacross}, is thus used for all ORAS solvers and preconditioners.  We set $\alpha=4$ and $\alpha^\times=40$ throughout.

%The need for an appropriate partitioning of the global domain when using ORAS preconditioners is exemplified by the sphere and torus problems when the subdomain count is increased.
Examining the results in Table \ref{tab:subAll}, we immediately see that the ORAS solvers and preconditioners outperform their RAS counterparts over the same splittings in all cases, often by a wide margin.
However, we see diminishing gains in performance for the ORAS
preconditioners over their RAS counterparts as the number of subdomains is increased.
For example,  for the sphere with $h=1/100$ and $N_S=32$ subdomains,
the iteration count is reduced from $53$ to $33$ by switching from RAS to ORAS.
The same comparison with $N_S=64$ subdomains gives a smaller relative reduction from $57$ iterations to $46$ iterations.
Using a large number of subdomains increases the number of cross points, and the modification then affects a larger fraction of the boundary nodes.
Raising the resolution of these problems to $h=1/200$ mitigates the issues with cross point density,
and the performance gap between the RAS and ORAS solvers widens for all surfaces and subdomain counts.

\begin{table}[htbp]
  \centering
  \begin{tabular}{@{}r|cc|cc|cc@{}}
    \toprule
    $N_S$            & \multicolumn{2}{c|}{32} & \multicolumn{2}{c|}{64} & \multicolumn{2}{c}{96}  \\
    Method           & RAS        & ORAS       & RAS        & ORAS       & RAS        & ORAS       \\ \midrule
    \multicolumn{7}{c}{Sphere, $h=1/100$, $N_A=654630$} \\
    Solver           &    500     &    245     &    704     &     278    &    818     &    405     \\
    Preconditioner   & 53         & 33         & 57         &  46        &  59        &  49        \\ \midrule
    \multicolumn{7}{c}{Sphere, $h=1/200$, $N_A=2618046$} \\
    Solver           &    916     &    209     &    1228    &     352    &    1318    &    434     \\
    Preconditioner   & 63         & 35         & 78         &  51        & 99         & 53         \\ \midrule
    \multicolumn{7}{c}{Torus, $h=1/100$, $N_A=454300$} \\
    Solver           & 597        & 192        & 780        & 549        & 1099       & 701        \\
    Preconditioner   & 53         & 33         & 57         & 49         & 63         & 55         \\ \midrule
    \multicolumn{7}{c}{Torus, $h=1/200$, $N_A=1818052$} \\
    Solver           & 1148       & 234        & 1511       & 343        & 2147       & 457        \\
    Preconditioner   & 67         & 38         & 83         & 48         & 104        & 57         \\ \midrule
    \multicolumn{7}{c}{Stanford Bunny, $h=1/100$, $N_A=396292$} \\
    Solver           & 717        & 344        & 1015       & 550        & 1306       & 771        \\
    Preconditioner   & 67         & 37         & 76         & 59         & 82         & 70         \\ \midrule
    \multicolumn{7}{c}{Stanford Bunny, $h=1/200$, $N_A=1585218$} \\
    Solver           & 1301       & 326        & 1906       & 528        & 2346       & 687        \\
    Preconditioner   & 80         & 38         & 140        & 65         & 143        & 73         \\ \bottomrule
  \end{tabular}
  \caption{The iterations to convergence as subdomain count varies for our three surfaces. Resolutions of $h=1/100$ and $h=1/200$ are presented. The ORAS solvers and preconditioners use $\alpha=4$ and $\alpha^\times=40$ throughout. The iterations to convergence consistently grow as $N_S$ is increased, as expected.}
  \label{tab:subAll}
\end{table}

All three surfaces show that the larger subdomain counts require more iterations.
The rise in iterations to convergence as the number of subdomains is increased is to be expected for any single-level method, such as those presented here. The use of a suitable coarse correction or multigrid method (cf. \cite{Chen:MG}) would likely mitigate this effect and would pose an interesting follow up to this work.

% -------------------------------------------------------------
\subsection{Parallel performance}
\label{results:parallel}
We conclude our evaluation of the (O)RAS methods with a study of the time to solution as the number of processors is increased.
All timings recorded here were run using whole nodes of the Graham compute cluster managed by Compute Canada.
This system has two 16 core Intel E5-2683 v4 Broadwell 2.1 GHz processors, and 125G of memory on each compute node.
Throughout this section, computations are carried out on the unit sphere with an overlap width of $N_O=4$; ORAS methods again use Robin parameters of $\alpha=4$ and $\alpha^\times=40$.

The timings presented are split into several phases:
\begin{itemize}
\item  global meshing, which includes building and partitioning the graph as well as scattering the permuted nodes across the processes,
\item  construction of the global matrix, which includes computation of the extension weights,
\item  construction and factorization of the local operators,
\item  solution of the system by the RAS and ORAS solvers, and
\item  solution of the system by GMRES preconditioned by (O)RAS.
\end{itemize}
The first three phases are performed before the solution is computed, and are thus independent of the number of iterations required by the solver or preconditioner.
The final two phases produce the same answer, so in practice only one would be needed. Each case was run $5$ times, with the average time of each phase reported below.

The global meshing procedures have not been fully parallelized in the current implementation of this software. The disjoint partitioning relies on METIS rather than ParMETIS, the parallel version of the software, which introduces an early serial bottleneck. As a consequence, our global meshing times do not show any parallel scalability. Note, however, that the serial components in the global meshing do not affect the subsequent steps which have been fully parallelized.

%% \begin{table}[htbp]
%%   \centering
%%   \begin{tabular}{rcccccc}
%%     \toprule
%%     & \multicolumn{3}{c}{RAS} & \multicolumn{3}{c}{ORAS} \\
%%     $N_{proc}$      & 32  & 64  & 128 & 32  & 64  & 128 \\ \midrule
%%     Global mesh     & 154 & 161 & 166 & 154 & 158 & 166 \\
%%     Global matrix   & 67  & 35  & 18  & 68  & 38  & 18  \\
%%     Local operators & 98  & 52  & 25  & 100 & 51  & 25  \\ \midrule
%%     Solver          & 612 & 385 & 192 & 165 & 93  & 46  \\
%%     Preconditioner  & 37  & 25  & 12  & 20  & 11  & 5   \\ \bottomrule
%%   \end{tabular}
%%   \caption{The time in seconds for each phase of the method for the sphere with grid spacing $h=1/250$ and  $N_S=128$ subdomains (independent of the number of processes). The first three phases are independent of the number of iterations used in the solution phase. Note  that only one of the last two phases needs to be used in practice. With the exception of the meshing procedure, all phases in both the RAS and ORAS methods show good strong scaling.}
%%   \label{tab:para128NoAl}
%% \end{table}

\begin{table}[htbp]
  \centering
  \begin{tabular}{rcccccc}
    \toprule
    & \multicolumn{3}{c}{RAS} & \multicolumn{3}{c}{ORAS} \\
    $N_{proc}$      & 32  & 64  & 128 & 32 & 64 & 128 \\ \midrule
    Global mesh     & 55  & 53  & 55  & 54 & 53 & 56  \\
    Global matrix   & 41  & 22  & 12  & 42 & 22 & 11  \\
    Local operators & 74  & 34  & 17  & 75 & 35 & 17  \\ \midrule
    Solver          & 348 & 276 & 196 & 49 & 68 & 37  \\
    Preconditioner  & 36  & 22  & 11  & 17 & 12 & 5   \\ \bottomrule
  \end{tabular}
  \caption{\revd{The time in seconds for each phase of the method for the sphere with grid spacing $h=1/250$ and  $N_S=N_{proc}$ subdomains. The first three phases are independent of the number of iterations used in the solution phase. Note that only one of the last two phases needs to be used in practice.}}
  \label{tab:paraStrNsp}
\end{table}

Strong scalability is quantified by measuring the change in run-time as the number of available processes is increased.  Ideally this run-time would decrease in proportion to the number of processes used. In our experiments, the subdomain count is set \revd{equal to the number of processes,} and the grid resolution is $h=1/250$, yielding a problem with $N_A=4098174$ active nodes. Aside from the global meshing phase, the method is easily parallelized and, as seen in Table~\ref{tab:paraStrNsp}, all \revd{set-up} phases of the method except for the meshing procedure show good strong scaling with the time spent in each phase roughly halving as the number of processors is doubled. \revd{The solver phase shows more complicated behavior. As the subdomain count is increased the solver performance becomes weaker, as noted in Section \ref{results:nsub}, and requires more iterations to converge. Simultaneously, the jump from $N_{proc}=32$ to $N_{proc}=64$ requires the use of two compute nodes and introduces communication latency. The preconditioned GMRES times do decrease with greater parallelism.} \revt{The results for the solver, in particular, highlight the need for a compatible coarse space and two level method.} The global meshing procedure requires a nearly constant time due to its serial nature.

\revt{The MUMPS \cite{MUMPS01,MUMPS02} direct solver was tested on the same problem used for the strong scalability test. When confined to a single compute node, where $N_{proc}=32$, the direct solver requires more memory than available, giving a run-time error. Increasing to two compute nodes, with $N_{proc}=64$, allows the direct solver to successfully run, but required $384$ seconds to find the factorization and solution. From Table \ref{tab:paraStrNsp} the time required for ORAS preconditioned GMRES on the same problem is only $35$ seconds for the assembly and factorization of the local operators, and $12$ seconds for GMRES to converge, for a total of $47$ seconds. The other phases of the run must be executed regardless of the choice of solution method. Similarly, with $N_{proc}=128$ the direct solver requires $325$ seconds for the relevant operations, while ORAS preconditioned GMRES requires just $23$ seconds. We conclude that iterative methods are preferable for problems of this scale.}

Weak scalability seeks to keep the load per processor nearly constant as $N_{proc}$ is increased, and ideally the time required would stay constant. The weak scalability of the method is evaluated by using one subdomain per process, $N_S=N_{proc}$, and using problem resolutions of $h=1/125, 1/175, 1/250$ for $N_{proc}=32, 64, 128,$ respectively. These resolutions have $N_A=1025622, 2005758, 4098174$ active nodes yielding approximately $32000$ nodes per subdomain. Table \ref{tab:para128Weak} records the time required by each phase of the method for this progression of problems. The construction of the global matrix and the local problems requires nearly constant time, indicating excellent weak scaling. However as $h$ is refined and $N_S$ is increased, the global system shows worse conditioning and the (O)RAS solver slows. These effects drive up the number of iterations required for convergence and disrupt the scalability of the solution phases. \revt{This is a characteristic of single-level solvers.}

\begin{table}[htbp]
  \centering
  \begin{tabular}{rcccccc}
    \toprule
    & \multicolumn{3}{c}{RAS} & \multicolumn{3}{c}{ORAS} \\
    $N_{proc}$       & 32 & 64 & 128 & 32  & 64 & 128 \\ \midrule
    Global mesh     & 12 & 25 & 55  & 13  & 25 & 56  \\
    Global matrix   & 11 & 11 & 12  & 10  & 11 & 11  \\
    Local operators & 17 & 17 & 17  & 16  & 17 & 17  \\ \midrule
    Solver          & 29 & 89 & 196 & 21  & 30 & 37  \\
    Preconditioner  & 5  & 7  & 11  & 4   & 5  & 5   \\ \bottomrule
  \end{tabular}
  \caption{The time in seconds for each phase of the method  for the sphere with varying grid spacings. The weak scalability of these methods is evaluated by choosing $N_S=N_{proc}$ and using $h=1/125, h=1/175, h=1/250$ for $N_{proc}=32, 64, 128,$ respectively. The first three phases are independent of the number of iterations used in the solution phase. Note that only one of the last two phases needs to be used in practice. The global meshing time grows in proportion to the problem size due to the sequential parts of that procedure.}
  \label{tab:para128Weak}
\end{table}

In summary, the methods exhibit good parallel scalability over the range of process counts tested. The shortcomings present in the weak scalability of the solution phases of the methods are to be expected for a single level method such as this. The meshing procedures are an obvious candidate for improvement, though this cost is fixed and would be amortized in an application requiring many solutions of the linear system.

% -------------------------------------------------------------
\section{Conclusion}
\label{conc}
This paper develops and evaluates RAS and ORAS domain decomposition solvers for the solution of surface intrinsic elliptic PDEs discretized by the CPM.
The splitting of the global problem and the construction of the local problems has been kept general through the use of a graph partitioner and the construction of effective boundary geometry that is independent of the irregular subdomains produced by the partitioning. To facilitate the implementation of ORAS methods, a new discretization of the Robin boundary operator for the CPM was given and verified to be first order accurate \revd{for a test case}. The presence of cross points in the splitting was identified as the source of instability of the ORAS iteration with small values of the Robin parameter. A modified scheme, using larger parameters in the vicinity of cross points, was presented and dramatically improves the performance of the ORAS methods. The solvers were shown to be a viable and effective approach to the parallelized solution of equations arising from a CPM discretization of surface intrinsic elliptic PDEs. \revt{As noted in Section \ref{results:parallel} RAS and ORAS are favored over sparse direct solvers, echoing the observations in \cite{Chen:MG}.}

Only single level methods were developed and evaluated herein, and \revt{the results of Section \ref{results:parallel} show the need for a two level solver to obtain good parallel scalability. The goal is to obtain convergence rates with reduced dependence on the number of subdomains. The addition of a coarse space compatible with the CPM is a particularly interesting direction for future work. For scalar elliptic and elasticity problems, much progress has been made on the introduction of coarse spaces; see \cite{dohrmann2017} for overlapping Schwarz methods and \cite{dohrmann2019} in the context of BDDC (balancing domain decomposition by constraints). For the non-overlapping optimized Schwarz case, see \cite{Dubois2012,harerssas2017}. Exploring these approaches in the embedding space is beyond the scope of this paper, but is part of ongoing work.}

\textit{Acknowledgments} The authors gratefully acknowledge the financial support of NSERC Canada (RGPIN 2016-04361 and RGPIN 2018-04881), and the preliminary work of Nathan King that helped inspire this project. This research was enabled in part by support provided by ACEnet (ace-net.ca) and Compute Canada.% (www.computecanada.ca).

\bibliographystyle{siamplain}
\bibliography{references}

\begin{thebibliography}{10}

\bibitem{MUMPS01}
{\sc P.~R. Amestoy, I.~S. Duff, J.-Y. L'Excellent, and J.~Koster}, {\em A fully
  asynchronous multifrontal solver using distributed dynamic scheduling}, SIAM
  Journal on Matrix Analysis and Applications, 23 (2001), pp.~15--41.

\bibitem{MUMPS02}
{\sc P.~R. Amestoy, A.~Guermouche, J.-Y. L'Excellent, and S.~Pralet}, {\em
  Hybrid scheduling for the parallel solution of linear systems}, Parallel
  Computing, 32 (2006), pp.~136--156.

\bibitem{petsc-efficient}
{\sc S.~Balay, W.~D. Gropp, L.~C. McInnes, and B.~F. Smith}, {\em Efficient
  management of parallelism in object oriented numerical software libraries},
  in Modern Software Tools in Scientific Computing, E.~Arge, A.~M. Bruaset, and
  H.~P. Langtangen, eds., Birkh{\"{a}}user Press, 1997, pp.~163--202.

\bibitem{Tref:Bary}
{\sc J.-P. Berrut and L.~N. Trefethen}, {\em Barycentric {Lagrange}
  interpolation}, SIAM Rev., 46 (2004), pp.~501--517.

\bibitem{Cheng:LSM}
{\sc M.~Bertalmio, L.-T. Cheng, S.~Osher, and G.~Sapiro}, {\em Variational
  problems and partial differential equations on implicit surfaces}, J. Comput.
  Phys., 174 (2001), pp.~759--780.

\bibitem{Cai:RAS}
{\sc X.-C. Cai and M.~Sarkis}, {\em A restricted additive {Schwarz}
  preconditioner for general sparse linear systems}, SIAM J. Sci. Comput., 21
  (1999), pp.~792--797.

\bibitem{Chen:MG}
{\sc Y.~Chen and C.~Macdonald}, {\em The closest point method and multigrid
  solvers for elliptic equations on surfaces}, SIAM J. Sci. Comput., 37 (2015),
  pp.~A134--A155.

\bibitem{Chu:Vari}
{\sc J.~Chu and R.~Tsai}, {\em Volumetric variational principles for a class of
  partial differential equations defined on surfaces and curves}, Res. Math.
  Sci., 5 (2018), p.~19.

\bibitem{dohrmann2019}
{\sc C.~R. Dohrmann, K.~H. Pierson, and O.~B. Widlund}, {\em Vertex-based
  preconditioners for the coarse problems of {BDDC}}, SIAM J. Sci. Comput., 41
  (2019), pp.~A3021--A3044, \url{https://doi.org/10.1137/19M1237557},
  \url{https://doi.org/10.1137/19M1237557}.

\bibitem{dohrmann2017}
{\sc C.~R. Dohrmann and O.~B. Widlund}, {\em On the design of small coarse
  spaces for domain decomposition algorithms}, SIAM J. Sci. Comput., 39 (2017),
  pp.~A1466--A1488, \url{https://doi.org/10.1137/17M1114272},
  \url{https://doi.org/10.1137/17M1114272}.

\bibitem{DoleanNataf}
{\sc V.~Dolean, P.~Jolivet, and F.~Nataf}, {\em An Introduction to Domain
  Decomposition Methods: Algorithms, Theory, and Parallel Implementation},
  SIAM, Philadelphia, PA, USA, 2015.

\bibitem{Dubois2012}
{\sc O.~Dubois, M.~J. Gander, S.~Loisel, A.~St-Cyr, and D.~B. Szyld}, {\em The
  optimized {S}chwarz method with a coarse grid correction}, SIAM J. Sci.
  Comput., 34 (2012), pp.~A421--A458, \url{https://doi.org/10.1137/090774434},
  \url{https://doi.org/10.1137/090774434}.

\bibitem{DziukElliot}
{\sc G.~Dziuk and C.~Elliott}, {\em Surface finite elements for parabolic
  equations}, Journal of Computational Mathematics, 25 (2007), pp.~385--407.

\bibitem{FloaterHormann:Para}
{\sc M.~S. Floater and K.~Hormann}, {\em Surface parameterization: a tutorial
  and survey}, Math. Vis. Advances in Multiresolution for Geometric Modelling,
  (2005), pp.~157--186.

\bibitem{Gand:OptPar}
{\sc M.~J. Gander}, {\em Optimized {Schwarz} methods}, SIAM J. Numer. Anal., 44
  (2006), pp.~699--731.

\bibitem{Gand:XPT}
{\sc M.~J. Gander and F.~Kwok}, {\em Best {Robin} parameters for optimized
  {Schwarz} methods at cross points}, SIAM J. Sci. Comput., 34 (2012),
  pp.~1849--1879.

\bibitem{harerssas2017}
{\sc R.~Haferssas, P.~Jolivet, and F.~Nataf}, {\em An adaptive coarse space for
  {P}. {L}. {L}ions algorithm and optimized {S}chwarz methods}, in Domain
  decomposition methods in science and engineering {XXIII}, vol.~116 of Lect.
  Notes Comput. Sci. Eng., Springer, Cham, 2017, pp.~43--53,
  \url{https://doi.org/10.1007/978-3-319-52389-7_4},
  \url{https://doi.org/10.1007/978-3-319-52389-7_4}.

\bibitem{METIS}
{\sc G.~Karypis and V.~Kumar}, {\em A fast and high quality multilevel scheme
  for partitioning irregular graphs}, SIAM J. Sci. Comput., 20 (1998),
  pp.~359--392.

\bibitem{Loisel:2Lag}
{\sc S.~Loisel}, {\em Condition number estimates for the nonoverlapping
  optimized {Schwarz} method and the 2-{Lagrange} multiplier method for general
  domains and cross points}, SIAM J. Numer. Anal., 51 (2013), pp.~3062--3083.

\bibitem{CBM:LSE}
{\sc C.~Macdonald and S.~Ruuth}, {\em Level set equations on surfaces via the
  closest point method}, Journal of Scientific Computing, 35 (2008),
  pp.~219--240.

\bibitem{CBM:Eig}
{\sc C.~B. Macdonald, J.~Brandman, and S.~J. Ruuth}, {\em Solving eigenvalue
  problems on curved surfaces using the closest point method}, J. Comput.
  Phys., 230 (2011), pp.~7944--7956.

\bibitem{CBM:RDonPC}
{\sc C.~B. Macdonald, B.~Merriman, and S.~J. Ruuth}, {\em Simple computation of
  reaction diffusion processes on point clouds}, Proceedings of the National
  Academy of Sciences, 110 (2013).

\bibitem{CBM:ICPM}
{\sc C.~B. Macdonald and S.~J. Ruuth}, {\em The implicit closest point method
  for the numerical solution of partial differential equations on surfaces},
  SIAM J. Sci. Comput., 31 (2010), pp.~4330--4350.

\bibitem{Maerz/Macdonald:cpfunctions}
{\sc T.~M\"{a}rz and C.~B. Macdonald}, {\em Calculus on surfaces with general
  closest point functions}, SIAM J. Numer. Anal., 50 (2012), pp.~3303--3328.

\bibitem{May:CPMCode}
{\sc I.~May}, {\em {DD-CPM}}.
\newblock \url{https://bitbucket.org/mayianm/dd-cpm/}, 2018.

\bibitem{May:DD25}
{\sc I.~{May}, R.~D. {Haynes}, and S.~J. {Ruuth}}, {\em {Domain Decomposition
  for the Closest Point Method}}, arXiv e-prints,  (2019), arXiv:1907.13606.
\newblock To appear in proceedings of 25th Domain Decomposition Meeting.

\bibitem{Argy:MoveSurf}
{\sc A.~Petras and S.~Ruuth}, {\em {PDEs} on moving surfaces via the closest
  point method and a modified grid based particle method}, J. Comput. Phys.,
  312 (2016), pp.~139--156.

\bibitem{QuarteroniValli}
{\sc A.~Quarteroni and A.~Valli}, {\em Domain decomposition methods for partial
  differential equations}, Numerical mathematics and scientific computation,
  Oxford University Press, 1999.

\bibitem{Reuter:ShapeDNA}
{\sc M.~Reuter, F.-E. Wolter, and N.~Peinecke}, {\em {Laplace--Beltrami}
  spectra as `{Shape-DNA}' of surfaces and solids}, Computer-Aided Design, 38
  (2006), pp.~342--366.

\bibitem{SJR:CPM}
{\sc S.~J. Ruuth and B.~Merriman}, {\em A simple embedding method for solving
  partial differential equations on surfaces}, J. Comput. Phys., 227 (2008),
  pp.~1943--1961.

\bibitem{Saad:ItMeth}
{\sc Y.~Saad}, {\em Iterative methods for sparse linear systems}, SIAM, 2009.

\bibitem{SmithBjorstadGropp}
{\sc B.~F. Smith, P.~E. Bj{\o}rstad, and W.~D. Gropp}, {\em Domain
  decomposition : parallel multilevel methods for elliptic partial differential
  equations}, Cambridge University Press, 1996.

\bibitem{Cyr:OMSORAS}
{\sc A.~St-Cyr, M.~J. Gander, and S.~J. Thomas}, {\em Optimized multiplicative,
  additive, and restricted additive {Schwarz} preconditioning}, SIAM J. Sci.
  Comput., 29 (2007), pp.~2402--2425.

\bibitem{ToselliWidlund}
{\sc A.~Toselli and O.~Widlund}, {\em Domain decomposition methods--algorithms
  and theory}, Springer Series in Computational Mathematics, 34, Springer,
  Berlin ; New York, 2005.

\bibitem{bunny}
{\sc G.~Turk and M.~Levoy}, {\em Zippered polygon meshes from range images}, in
  Proceedings of the 21st annual conference on computer graphics and
  interactive techniques, SIGGRAPH '94, ACM, July 1994, pp.~311--318.

\end{thebibliography}

\end{document}